\documentclass[a4paper,11pt]{article}

\usepackage{amsmath, amssymb, amsfonts, xypic , amsthm, amscd, enumerate, vmargin}
\usepackage[utf8]{inputenc}
\usepackage[french]{babel}

\setmarginsrb{2.4cm}{3.5cm}{2.4cm}{3.3cm}{0cm}{0cm}{0cm}{2.5cm}

\theoremstyle{plain}
\newtheorem{teo}{Th\'eor\`eme}[section]
\newtheorem{lem}[teo]{Lemme}
\newtheorem{prop}[teo]{Proposition}
\newtheorem{cor}[teo]{Corollaire}

\theoremstyle{definition}

\newtheorem{rem}[teo]{Remarque}
\newtheorem{rems}[teo]{Remarques}

\newcommand{\Spec}{\operatorname{Spec}}

\newcommand{\Sym}{\operatorname{Sym}}
\newcommand{\Tr}{\operatorname{Tr}}
\newcommand{\tr}{\operatorname{tr}}

\renewcommand{\det}{\operatorname{det}}

\newcommand{\Cbb}{{\mathbb C}}

\newcommand{\Zbb}{{\mathbb Z}}
\newcommand{\Nbb}{{\mathbb N}}
\newcommand{\Ocal}{{\mathcal O}}

\newcommand{\Lcal}{{\mathcal L}}
\newcommand{\Ncal}{{\mathcal N}}
\newcommand{\Tcal}{{\mathcal T}}
\newcommand{\Vcal}{{\mathcal V}}
\newcommand{\Pbb}{{\mathbb P}}
\newcommand{\Abb}{{\mathbb A}}
\newcommand{\Rbb}{{\mathbb R}}
\newcommand{\Hbb}{{\mathbb H}}
\newcommand{\Obb}{{\mathbb O}}
\newcommand{\Mcal}{{\mathcal M}}
\newcommand{\Gbb}{{\mathbb G}}
\newcommand{\Fcal}{{\mathcal F}}

\author{Sébastien Jansou}
\title{Déformations des cônes de vecteurs primitifs}
\begin{document}
\date{}
\maketitle

\begin{abstract}
~\\
Pour un groupe réductif connexe complexe $G$, on classifie les modules simples dont le cône des vecteurs primitifs admet une déformation $G$-invariante non triviale. On relie cette classification à celle des algèbres de Jordan simples, et aussi à celle (due à Akhiezer) des variétés projectives lisses dont les orbites sous l'action d'un groupe algébrique affine connexe sont un diviseur et son complémentaire.

Notre principal outil est le schéma de Hilbert invariant d'Alexeev et Brion; on en détermine les premiers exemples.

On détermine aussi les déformations infinitésimales (non nécessairement $G$-invariantes) des cônes des vecteurs primitifs; elles sont triviales pour presque tous les modules simples.\\

\begin{center} \textbf{Abstract}
\end{center}

For a complex connected reductive group $G$, we classify the simple modules whose cone of primitive vectors admits a nontrivial $G$-invariant deformation. We relate this classification to that of simple Jordan algebras, and to that (due to Akhiezer) of smooth projective varieties whose orbits under the action of a connected affine algebraic group are a divisor and its complementary.

Our main tool is the invariant Hilbert scheme of Alexeev-Brion; we determine the first examples of it.

We also determine the infinitesimal deformations (non necessarily $G$-invariant) of the cones of primitive vectors; they turn out to be trivial for most simple modules.

\end{abstract}

\section*{Introduction}

Soit $G$ un groupe réductif connexe complexe, et $V$ un $G$-module rationnel de dimension finie. On dit qu'un sous-schéma fermé $G$-stable $X \subseteq V$  est \textit{à multiplicités finies} si son algèbre affine est somme directe de modules simples avec des multiplicités finies.\\

Alexeev et Brion ont montré récemment dans \cite{AB} que les sous-schémas $X$ ayant des multiplicités finies fixées sont paramétrés par un schéma quasi-projectif: \textit{le schéma de Hilbert invariant}. Leur travail est basé sur celui de Haiman et Sturmfels (\cite{HaS}), qui correspond au cas particulier où le groupe réductif $G$ est un tore.\\

On se propose ici de déterminer le schéma de Hilbert invariant dans le cas ``le plus simple''. L'espace ambiant est le $G$-module simple de plus grand poids $\lambda$:$$V=V(\lambda).$$ On cherche à paramétrer les plus petits sous-schémas fermés $G$-stables de $ V$ de dimension positive: on va donc choisir les multiplicités les plus petites possibles. Pour cela, on remarque que l'algèbre affine d'un tel schéma $X$ contient le dual de $V(\lambda)$; on note $\lambda^*$ son plus grand poids. Notons $f$ un vecteur de plus grand poids de $V(\lambda ^ *)$. Les puissances de $f$ sont non nulles dans l'algèbre affine de $X$, donc celle-ci contient tous les modules simples $V(d \lambda ^ *)$, où $d$ est un entier positif.\\

Ainsi, on va prendre pour multiplicité $1$ pour les modules simples $V(d \lambda ^ *)$, et $0$ pour les autres modules simples. Un point particulier du schéma de Hilbert invariant correspondant $H_\lambda$ est alors donné par le cône des vecteurs primitifs de $V(\lambda)$ (réunion de l'orbite des vecteurs de plus grand poids, et de l'origine). 
Ces cônes ne sont autres que les cônes affines sur les variétés de drapeaux plongées par un système linéaire complet.\\

On montrera que pour la plupart des poids $\lambda$, le schéma $H_\lambda$  est en fait réduit à ce point. Dans les autres cas, $H_\lambda$ est la droite affine. On obtient ainsi une classification des $G$-modules simples dont le cône des vecteurs primitifs admet une déformation invariante.\\

Cette classification rappelle celle (obtenue par Akhiezer dans \cite{Ak1}) des variétés projectives lisses dont les orbites sous l'action d'un groupe algébrique affine connexe sont un diviseur ample et son complémentaire. On décrira comment relier ces deux classifications et on les reliera aussi à celle des algèbres de Jordan simples. Cette dernière se trouve être plus ``petite'', mais on verra que lorsqu'un cône de vecteurs primitifs admet une déformation invariante dans $V(\lambda)$, cette déformation provient d'une algèbre de Jordan simple.\\

On montrera enfin que pour la plupart des poids $\lambda$, le cône des
vecteurs primitifs de $V(\lambda)$ est en fait rigide: les seuls cônes de
vecteurs primitifs admettant des déformations infinitésimales non triviales
sont, outre ceux qui admettent une déformation invariante, le cône affine
$C_m$ au
dessus de la courbe rationnelle normale de degré $m$ dans $\Pbb ^m$ (les
déformations de ce cône ont été étudiées dans \cite{Pink}), et le cône affine
$C_{m n}$
au dessus de $\Pbb ^1 \times \Pbb ^n$ dans le plongement de bidegré $(m,1)$
avec $m \geq 2$. On verra que la déformation verselle de $C_{m n}$ est lisse,
contrairement à celle de $C_m$.\\

\noindent \textit{Remerciements.} Je suis très vivement reconnaissant envers mon directeur de thèse Michel Brion pour son aide tout au long de ce travail.

\section{Une classe de schémas de Hilbert invariants }

\subsection{Notations et résultat principal}\

On considère des schémas et des groupes algébriques sur $\Cbb$.
Les références utilisées sont \cite{H} pour la théorie des schémas et \cite{PV} pour celle des groupes algébriques de transformations.

Soit $G$ un groupe réductif connexe. On en choisit un sous-groupe de Borel $B$, et un tore maximal $T$ inclus dans $B$. On considère le radical unipotent $U$ de $B$ : on a : $B=TU$.
Les algèbres de Lie respectives de $G$, $B$, $T$ et $U$ sont notées : $\frak{g}$, $\frak{b}$, $\frak{t}$, et $\frak{u} $.
Le système de racines de $G$ relativement à $T$ est noté $R$.
Le choix de $B$ nous en fournit une base $S$, et on a : $R=R_+ \amalg R_-$ où $R_+$ est l'ensemble des racines positives, et $R_-$ celui des racines négatives.

On note $\Lambda$ le groupe des caractères de $T$. On a un ordre partiel sur $\Lambda$ : $\mu \leq \lambda$ si et seulement si $\lambda - \mu$ est une somme de racines positives.
On note $\Lambda^+$ l'ensemble des éléments de $\Lambda$ qui sont des poids dominants (relativement à la base $S$ du système de racines $R$).
On sait que $\Lambda ^ +$ est en bijection avec l'ensemble des classes d'isomorphisme de $G$-modules rationnels simples. Si $\lambda$ est un élément de $\Lambda^+$, on notera $V(\lambda)$ un $G$-module simple correspondant, c'est-à-dire de plus grand poids $\lambda$, et $v_\lambda$ un vecteur de $V( \lambda)$ de poids $\lambda$. Si $\lambda$ est un poids qui n'est pas dominant, on pose $V(\lambda) =0$. 

Les $G$-modules simples peuvent être construits de la fa{\c c}on suivante : soit $\lambda$ un poids dominant, et $P$ un sous-groupe parabolique de $G$ contenant $B$ tel que $\lambda$ se prolonge en un caractère de $P$. Notons $\pi : G \rightarrow G/P$ la surjection canonique, et $\Lcal _ \lambda$ le faisceau inversible sur $G/P$ qui associe à un ouvert $\Omega \subseteq G/P$ : $$\Lcal _ \lambda ( \Omega ) := \{ f \in \Ocal _G ( \pi ^ {- 1} ( \Omega) ) \mbox{ $\vert$ } \forall g \in G \mbox{, } \forall p \in P \mbox{, } f(gp)= \lambda ( p ) f (g) \}.$$
Le faisceau $\Lcal _ \lambda$ est alors $G$-linéarisé (via l'action de $G$ sur ses fonctions régulières par translation à gauche), et l'espace des sections globales de $\Lcal _ \lambda$ est un $G$-module simple :$$ \Gamma ( G/P , \Lcal _ \lambda) \simeq V ( \lambda )^* . $$

Si $V$ est un $T$-module rationnel (éventuellement de dimension infinie), on note $V=\bigoplus_{ \lambda \in  \Lambda} V_\lambda $ sa décomposition en sous-espaces propres.
Par exemple, l'algèbre de Lie de $G$ admet la décomposition : $$\frak{g}=\displaystyle{  \frak{t}  \oplus \bigoplus_{ \alpha \in R}\frak{g}_\alpha },$$ où chaque $\frak{g}_\alpha$ est de dimension 1. On choisit pour tout  $ \alpha \in R $ un générateur $e_\alpha$ de $\frak{g}_\alpha$.

Si $V$ est un $G$-module rationnel, on note  $V_{(\lambda)}$ sa composante isotypique de type $ \lambda$, c'est-à-dire le sous-module de $V$ somme des sous-modules isomorphes à $ V(\lambda ) $. On a alors la décomposition $V=\bigoplus_{ \lambda \in  \Lambda ^ +} V_{(\lambda)} $.

Soit $V$ un $G$-module rationnel de dimension finie, et $h :\Lambda ^ + \longrightarrow \Nbb$ une fonction.
On appelle \textit{famille de sous-schémas fermés $G$-stables de $V$} un sous-schéma fermé $G$-stable de  $ \frak{X} \subseteq S \times V$, où $S$ est un schéma avec action triviale de $G$. On note $\pi  : \frak{X} \rightarrow S $ le morphisme induit par la projection $S \times V \rightarrow S$, et $\pi _ * \Ocal _ { \frak{X} }$ l'image directe par $\pi$ du faisceau structural de $\frak{X}$ .
La famille $\frak{X}$ est dite \textit{de fonction de Hilbert $h$} si on a un isomorphisme de $\Ocal _ S$- $G$-modules $$ \pi _ * \Ocal _ { \frak{X} } \simeq \bigoplus _ {\lambda \in \Lambda^+} \Fcal _ \lambda \otimes V ( \lambda),$$ 
 où chaque $\Fcal _ \lambda$ est un $\Ocal _ S$-module localement libre de rang $h(\lambda)$. (Le morphisme $\pi$ est alors plat.)

Le foncteur contravariant : $( \mbox{Schémas}) ^ \circ \longrightarrow (\mbox{Ensembles})$ qui associe à tout schéma $S$ l'ensemble des familles $ \frak{X} \subseteq  S \times V$ de fonction de Hilbert $h$ est représenté par un schéma quasi-projectif noté $Hilb_h^G(V)$. (On renvoie à \cite{AB}\S1.2 pour plus de détails.)\\

On fixe désormais un poids dominant $\lambda$.
On note $\lambda^*$ le plus grand poids du $G$-module $V(\lambda) ^ *$ dual de $V(\lambda)$.
Soit $h_\lambda :\Lambda^+ \rightarrow \Nbb$  la fonction valant 1 sur $\Nbb \lambda^*$ et $0$ ailleurs.
On note dans la suite $$H_\lambda  := Hilb_{h_\lambda}^G(V(\lambda))$$ le schéma de Hilbert invariant associé à ce choix.

Si $E$ est un espace vectoriel de dimension finie, on note $\Pbb ( E ) $ l'espace de ses droites.
On a une action régulière de $G$ sur l'espace $\Pbb (V( \lambda)  ) $.
Notons $[ v _ \lambda ] \in \Pbb (V( \lambda)  ) $ la droite engendrée par $v _ \lambda$ et $$P _ \lambda  := G_{[ v _ \lambda ]}$$son stabilisateur dans $G$ : c'est le plus grand sous-groupe parabolique de $G$ qui contient $B$ et tel que $\lambda$ se prolonge en un caractère de $P_\lambda$.
L'orbite de $[ v _ \lambda ]$ est la seule orbite fermée de $\Pbb (V( \lambda)  ) $ (donc l'unique orbite de plus petite dimension).
L'espace homogène projectif $G / P_ \lambda $ se plonge ainsi dans $\Pbb (V( \lambda)  ) $, et le faisceau inversible très ample associé à ce plongement est en fait $\Lcal _ \lambda$.
Le cône affine au dessus de $G / P_ \lambda $ dans $V( \lambda)$ est le cône $$C_\lambda := G.v_\lambda \cup \{ 0 \}   =\overline{G.v_\lambda}$$des vecteurs primitifs de $V( \lambda)$. C'est une variété normale (\textit{cf.} \cite{Kra}, III.3.5). La variété $G / P_ \lambda \subseteq \Pbb (V( \lambda))$ est donc projectivement normale, et l'algèbre affine graduée du cône $C_\lambda$ est $$ \bigoplus_{d \in \Nbb} \Gamma (G / P_ \lambda, \Lcal _ {d \lambda }) =  \bigoplus_{d \in \Nbb}  V( d \lambda^*).$$ On peut donc voir $C_\lambda$ comme un point fermé de $H_\lambda$.

L'objet de cette partie est de montrer le théorème suivant, énoncé avec les notations de \cite{Bou} :

\begin{teo}
Le schéma de Hilbert invariant $ H_\lambda$ est un point réduit, sauf dans les cas suivants où $ H_\lambda$ est la droite affine :

(H1) $G$ est simple de type $A_1$, et $\lambda = 2 \omega_1 \mbox{ ou } 4 \omega_1$.

(H2) $G$ est simple de type $A_n$, $n \geq 2$ et $\lambda =  \omega_1 + \omega_n$. 

(H3) $G$ est simple de type $B_3$ et $\lambda = \omega_3 \mbox{ ou } 2 \omega_3$.

(H4) $G$ est simple de type $B_n$, $n \geq 2$ et $\lambda = \omega_1 \mbox{ ou } 2 \omega_1$.

(H5) $G$ est simple de type $C_n$, $n \geq 3$ et $\lambda = \omega_2$.

(H6) $G$ est simple de type $D_n$, $n \geq 3$ et $\lambda = \omega_1 \mbox{ ou } 2 \omega_1$.

(H7) $G$ est simple de type $F_4$ et $\lambda = \omega_4 $.

(H8) $G$ est simple de type $G_2$ et $\lambda = \omega_1 \mbox{ ou } 2 \omega_1$.

(H9) $G$ est semi-simple de type $A_1 \times A_1 $ et $\lambda = ( \omega_1,\omega_1)\mbox{ ou }  ( 2 \omega_1, 2 \omega_1)  $.

\noindent et dans les cas $(G, V(\lambda))$ obtenus à partir d'un cas $(G_0, V_0)$ parmi les précédents par factorisation : $G \rightarrow G_0 \rightarrow GL( V_0)$.
\end{teo}

\begin{rem}\
\begin{itemize}
 \item Le cas (H6) avec $n=3$ revient à un groupe $G$ simple de type $A_3$, avec $\lambda = \omega_2 \mbox{ ou } 2 \omega_2$.
 \item On pourrait voir le cas (H9) comme étant le cas (H6) avec $n=2$.
 \end{itemize}
\end{rem}

\subsection{Action du groupe multiplicatif sur le schéma de Hilbert invariant }\

On a une opération naturelle du groupe multiplicatif  $\Gbb_m$ sur le schéma de Hilbert invariant $H_\lambda$ : elle provient de l'action de $\Gbb_m$ sur $V(\lambda)$ par homothéties (qui commute avec l'action de $G$).

Dans cette partie, on montre que cette action admet pour unique point fixe le cône $C_\lambda$ des vecteurs primitifs. On montre aussi que $C_\lambda$ est dans l'adhérence de toutes les orbites de $\Gbb_m$, et on en déduit que le schéma de Hilbert invariant est affine.
Ces propriétés peuvent être déduites de ce qui est fait dans \cite{AB}, \S 2.1 à 2.3~; on a préféré donner ici des preuves directes.

\begin{prop}

(a) Le cône $C_\lambda$ est l'unique point fermé de $H_\lambda$ fixé par $\Gbb_m$.\\
(b) Soit $X$ un point fermé de $H_\lambda$.
    Le morphisme : $\Gbb_m \longrightarrow H_\lambda$, $t \longmapsto t.X$ se prolonge en un morphisme $\Abb ^1 \longrightarrow H_\lambda$, $0 \longmapsto C_\lambda$. 
\end{prop}

\noindent {\bf Preuve.} (a) On note $S^e  V( \lambda)^*$ la puissance symétrique d'ordre $e$ de $V( \lambda)^*$.
On identifie l'algèbre des fonctions régulières sur $V( \lambda)$ à l'algèbre symétrique de $V( \lambda)^*$ : $$ \Sym V( \lambda)^*  := \bigoplus _ {e \in \Nbb} S^e  V( \lambda)^* .$$ 
Les points fermés fixés par $\Gbb_m$ correspondent aux idéaux homogènes $$I= \bigoplus _{e \in \Nbb} I_e \subseteq  \Sym V(\lambda)^* = \bigoplus _{e \in \Nbb} S^e V(\lambda)^*$$ qui sont stables par $G$ et de fonction de Hilbert $h$.

On sait que $S^e V(\lambda)^*$ contient un unique sous-$G$-module isomorphe à $V( e \lambda)^*$, et que ses autres composantes isotypiques non nulles sont de type inférieur à $e \lambda ^* $ :
\begin{equation}\label{2}
 S^e V(\lambda)^* \simeq V(e \lambda)^* \oplus \bigoplus _ { \mu < e \lambda ^ *}[ S^e V(\lambda)^*]_{( \mu)}.      
\end{equation}
On va montrer par récurrence sur $e \in \Nbb $, qu'un tel idéal $I$ vérifie :
 
\begin{equation}\label{3}
I_e= \bigoplus _ { \mu < e \lambda ^ *}[ S^e V(\lambda)^*]_{( \mu)}.   
\end{equation}
En effet, comme $I$ ne contient pas les constantes, on a $I_0=0$.
Puis, si $(2)$ est satisfait pour tout $d<e$, il faut que :
 $$  \bigoplus _ { \mu < e \lambda ^ *}[ S^e V(\lambda)^*]_{( \mu)} \subseteq  I_e$$ pour que la fonction de Hilbert de $I$ soit $h_\lambda$.
Cette dernière inclusion est en fait une égalité, car sinon on aurait $I_d = S^d V(\lambda)^*$ pour tout $d\geq  e$.

Il n'y a donc pas d'autre point fermé fixé par $\Gbb_m$ que $C_\lambda$.

(b) Pour mieux comprendre l'action de $\Gbb_m$ sur $H _ \lambda$, on reprend la construction du schéma de Hilbert invariant (voir \cite{HaS}, \S 1,2,3 et \cite{AB}, \S 1.2).

Considérons l'action naturelle de  $\Gbb_m$ sur l'algèbre symétrique de $V( \lambda)^*$, où $\Gbb_m$ opère sur la composante  $S^e  V( \lambda)^*$ avec le poids $-e$, de sorte que $ \Sym V( \lambda)^*$ est une $G \times  \Gbb_m $-algèbre rationnelle.

La sous-algèbre $ [ \Sym V( \lambda)^* ]^U $ des invariants par $U$ est alors une $T \times \Gbb_m$- algèbre rationnelle de type fini, selon \cite{Gross},Thm 9.4.

On en choisit un système fini de générateurs $ f_1,...f_n$ formé de $T \times \Gbb_m$-vecteurs propres, et on note $S= \Cbb [ x_1,...,x_n]$ l'algèbre de polynômes correspondante. L'algèbre $S$ est naturellement une $T \times \Gbb_m$- algèbre rationnelle, et on a un morphisme surjectif de  $T \times \Gbb_m$- algèbres rationnelles :
$$\pi : S \longrightarrow ( \Sym V( \lambda)^*)^U.$$

L'action de $T$ sur $S$ fournit une graduation de $S$ par le groupe abélien $\Lambda$. 

On peut alors identifier $H_\lambda$ à un sous-schéma localement fermé d'un produit de Grassmanniennes, donc d'un produit d'espaces projectifs.
Plus précisément, on sait (\cite{HaS}) qu'il existe une partie finie $D$ de  $\Lambda$, et pour tout $\mu \in D$, un sous-espace vectoriel de dimension finie $N_\mu$ de $S_\mu$ que l'on peut choisir stable par $\Gbb_m $, tels que l'on ait un plongement :
\begin{equation}\label{1}
 H_\lambda  \hookrightarrow  \prod _{\mu \in D} \Pbb ({ \bigwedge} ^ {r_\mu} N_\mu)
\end{equation} où $ r_\mu  := \dim N_\mu - h_\lambda (\mu)$.
Décrivons l'image d'un point fermé par ce plongement : si $I$ est l'idéal d'un sous-schéma fermé de $V( \lambda )$ correspondant à un point fermé de $ H_\lambda$, on lui associe pour tout $\mu \in D$ : 
$$ J_\mu  := \pi ^ {-1} ( I ^ U) \cap N _\mu . $$

Les $N _\mu$ sont des modules rationnels pour l'action de $\Gbb_m$, donc les  $ \Pbb (\wedge ^ {r_\mu} N _\mu )$ sont munis d'une action régulière de $\Gbb_m$, pour laquelle le plongement $(3)$ est équivariant. 

On peut maintenant vérifier le point (b) de la proposition :\

Soit $\mu \in D$. Si $\mu \notin \Nbb \lambda ^ *$, alors $ r_\mu = \dim N_\mu$, et $ \Pbb (\wedge ^ {r_\mu} N _\mu )$ est réduit à un point.
Sinon, écrivons  $ \mu = e \lambda ^*$. On a $ \Pbb (\wedge ^ {r_\mu} N _\mu ) \cong  \Pbb ( N _\mu ^* )$. Notons $K  := N_\mu \cap \ker \pi$ et $L$ un supplémentaire $\Gbb_m $-stable de $K$ dans $N_\mu$ :$$N_\mu = K \oplus L.$$ Selon la décomposition (1) (considérée à tous les ordres), le plus grand poids de l'action de $\Gbb _ m$ sur $L$ est $-e$ :$$ L= L _ {-e} \oplus \bigoplus _{c>e} L_{-c}.$$ L'espace vectoriel $J_\mu$ est un hyperplan de $K \oplus L$ qui contient $K$ (par définition de $J_\mu$ et $K$) mais qui ne contient pas $L_{-e} $ selon le lemme qui suit.

Montrons alors que le morphisme $\Gbb_m \longrightarrow  \Pbb ( N_\mu^*)$, $t \longmapsto t.J_\mu$ se prolonge en un morphisme $f : \Abb ^1 \longrightarrow  \Pbb ( N_\mu^*)$ en posant $f(0) :=  K \oplus \bigoplus _{c>e} L_{-c}$.
Choisissons une base de $N _\mu ^*$ compatible avec la décomposition $N_\mu= K \oplus  L _ {-e} \oplus \bigoplus _{c>e} L_{-c}$. Notons $d$ la dimension de $L$.
Les coordonnées homogènes de $J_\mu$ dans $\Pbb ( N _\mu ^* )$ sont $[\underbrace{0 :... :0}_{ \mbox{\scriptsize{$\operatorname{dim}(K)$ fois}}} :x_1 :x_2 :... :x_d]$ et on a $x_1 \neq 0$. Celles de $t.J_\mu$ sont donc $[0 :... :0 :t^e x_1 :t^{c_2} x_2 :... :t^{c_d} x_d]$, où les $c_j$ sont des entiers strictement supérieurs à $e$. D'où l'assertion.

Comme $f(0)$ ne dépend pas de l'idéal $I$ considéré, il s'agit de $\pi ^ {-1} ( I_0 ^ U) \cap N _\mu$ où $I_0$ est l'idéal du cône $C_\lambda$, d'où (b).
\hfill $\Box$

\begin{lem}

Soit $X$ un sous-schéma fermé de $V( \lambda )$ de fonction de Hilbert 
$h_\lambda$, et $I \subseteq \Sym V( \lambda ) ^*$ son idéal.\\
Alors pour tout $e \in \Nbb$, le sous-$G$-module de $S^e V(\lambda)^*$ isomorphe à $V(e\lambda)^*$ n'est pas inclus dans $I$.
\end{lem}

\noindent {\bf Preuve.} Notons $f$ un vecteur de plus grand poids de $V(\lambda)^*$.

Le $G$-module $ [S^e V(\lambda)^*]_{(e \lambda ^ *)}$ est simple, et $f^e$ en est un vecteur de plus grand poids.
Supposons par l'absurde : $f^e \in I$.
Alors $f$ appartient à l'idéal du sous-schéma réduit $X_{red}$ associé à $X$.
Le sous-espace vectoriel de $ V(\lambda)$ engendré par $X_{red}$ est un sous-$G$-module de $ V(\lambda)$ inclus dans l'hyperplan défini par $f$ : il est donc réduit à $ \{ 0 \} $, et l'espace vectoriel $\Cbb [X]$ est de dimension finie : une contradiction.
\hfill $\Box$

\begin{cor}
Le schéma $H_\lambda$ est affine. Son algèbre affine $A$ est graduée par l'action du groupe multiplicatif sur $H_\lambda$, en degrés négatifs :  $A = \bigoplus _ {d \in - \Nbb} A_d$, et l'anneau $A_0$ est local.
\end{cor}

\noindent {\bf Preuve.}
On rappelle que $H_\lambda$ s'identifie à un sous-schéma localement fermé  $\Gbb_m$-stable de $ \Pbb ( M)$, où $M$ est un $\Gbb_m$-module rationnel de dimension finie.
Le sous-schéma réduit $\bar{H_\lambda}\setminus H_\lambda$ est donc aussi $\Gbb_m$-stable, et son idéal homogène aussi.
Il existe donc un élément homogène $f \in \Sym ( M ^ * )$, $\Gbb_m$-vecteur propre, définissant un ouvert $U$ contenant $C_\lambda$.
Cet ouvert $U$ est $\Gbb_m$-stable, et contient donc $H_\lambda$, selon le point (b) de la proposition précédente.
Ainsi, $H_\lambda$ est fermé dans $U$, et $U$ est un ouvert affine : $H_\lambda$ est affine.

Montrons maintenant que $ A_e = 0$ pour tout $e\geq 0$.
Par l'absurde, soit $f \in A \setminus \{ 0 \}$ de degré $e>0$.
Soit $X$ un point de $ H_\lambda$ tel que $f(X) \not = 0$.
La fonction $\Gbb_m \longrightarrow \Cbb$, $t \longmapsto f(tX)=t^{-e}f(X)$ se prolonge en une fonction régulière sur $\Abb ^ 1$ : une contradiction.

Enfin, montrons que  $A_0=A^{\Gbb_m }$ n'a qu'un seul idéal maximal.
On sait (\cite{Kra}, II.3.2) que le morphisme $\Spec (A) \rightarrow \Spec (A^{\Gbb_m })$ est surjectif. Donc si $A^{\Gbb_m }$ avait deux points fermés distincts, $H_\lambda$ aurait deux fermés $\Gbb_m$-stables disjoints : une contradiction avec le point (b) de la proposition précédente.
\hfill $\Box$

\subsection{Calcul de l'espace tangent au point fixe}

L'objet de ce paragraphe est de montrer le résultat suivant :

\begin{prop}
L'espace tangent $T_{C_\lambda} H_\lambda $ est nul, sauf dans les cas (H1) à (H9) du théorème 1.1 où il est de dimension 1.

\end{prop}

Notons $G_{v_\lambda}$ le stabilisateur de $v_\lambda$ dans $G$. L'espace tangent en $v_ \lambda$ à $C_\lambda= G.{v_\lambda} \cup \{ 0 \}$ est $\frak{g}.v_\lambda \subseteq V(\lambda)$. Il est stabilisé par l'action de $G _ {v_\lambda}$. On note enfin $[V(\lambda)/ \frak{g}.v_\lambda]^{G_{v_\lambda}}$ l'espace des invariants du quotient par $G_{v_\lambda}$.

Le point de départ de la démonstration est l'isomorphisme canonique :$$T_{C_\lambda} H_\lambda \cong [V(\lambda)/ \frak{g}.v_\lambda]^{G_{v_\lambda}}.$$
Cet isomorphisme découle de la proposition 1.5 (iii) de \cite{AB}. En effet on peut supposer que l'espace vectoriel $V(\lambda)$ n'est pas une droite : la variété $C_\lambda$ est alors de dimension supérieure ou égale à $2$, et normale. La codimension de $C_\lambda \backslash G.v_\lambda = \{ 0 \}$ est donc supérieure ou égale à $2$, et la proposition s'applique.\\

\begin{lem}
On a $G_{v_\lambda} = T_{v_\lambda} . G_{v_\lambda}^\circ $,
en notant $T_{v_\lambda}$ le stabilisateur de $v_\lambda$ dans $T$ (on a $T_{v_\lambda}= \operatorname{ker}(\lambda)$ et $G_{v_\lambda}^\circ $ la composante neutre de $G_{v_\lambda}$.
\end{lem}

\noindent {\bf Preuve.} Considèrons la décomposition de Lévi de $P_\lambda = G_{[v_\lambda]}$  relative à $T$ : $$P_\lambda = L_\lambda.U_\lambda.$$
Comme $T$ est un tore maximal du groupe réductif $L_\lambda$, on a $ L_\lambda = T . [L_\lambda,L_\lambda]$.
D'où $P_\lambda = T .  [L_\lambda,L_\lambda]. U_\lambda$
et  $G_{v_\lambda} = T_{ v_\lambda }.  [L_\lambda,L_\lambda]. U_\lambda$ (car $ [L_\lambda,L_\lambda] $ et $U_\lambda$ stabilisent $v_\lambda$), d'où le résultat, car $ [L_\lambda,L_\lambda] $ et $U_\lambda$ sont connexes.
\hfill $\Box$

\begin{prop}

On a une action du tore $T$ sur l'espace $ [V(\lambda)/ \frak{g}.v_\lambda]^{G_{v_\lambda}} $. Sa décomposition en sous-espaces propres est : $$[V(\lambda)/ \frak{g}.v_\lambda]^{G_{v_\lambda}} =
[V(\lambda)/ \frak{g}.v_\lambda]^{U}_0 \oplus
[V(\lambda)/ \frak{g}.v_\lambda]^{U}_{- \lambda}$$

\end{prop}

\noindent {\bf Preuve.} On observe que les poids de $V(\lambda)$ qui sont des multiples de $\lambda$ sont $\lambda$ et éventuellement $0$ et $-\lambda$.

Comme $P_\lambda$ stabilise  $\frak{g}.v_\lambda $, il agit sur $ V(\lambda)/ \frak{g}.v_\lambda$ et sur  $ [V(\lambda)/ \frak{g}.v_\lambda]^{G_{v_\lambda}}$ et  $ [V(\lambda)/ \frak{g}.v_\lambda]^{G_{v_\lambda}^\circ}$ (puisque $G_{v_\lambda}$ et  $ G_{v_\lambda}^\circ$ sont des sous-groupes distingués de $ P_\lambda$).

On a donc une action du tore $T$ sur $ [V(\lambda)/ \frak{g}.v_\lambda]^{G_{v_\lambda}}$, et ses poids sont des poids de $V( \lambda ) $ qui sont multiples de $\lambda$, car la restriction de l'action à $T_{v_\lambda}$ est triviale : $$[V(\lambda)/ \frak{g}.v_\lambda]^{G_{v_\lambda}} =
[V(\lambda)/ \frak{g}.v_\lambda]^{G_{v_\lambda}}_0 \oplus
[V(\lambda)/ \frak{g}.v_\lambda]^{G_{v_\lambda}}_{- \lambda}$$
(le poids $\lambda$ n'apparaît pas car $V(\lambda)_\lambda$ est inclus dans $ \frak{g}.v_\lambda $). D'où, selon le lemme précédent $$[V(\lambda)/ \frak{g}.v_\lambda]^{G_{v_\lambda}} =
[V(\lambda)/ \frak{g}.v_\lambda]^{G_{v_\lambda}^\circ}_0 \oplus
[V(\lambda)/ \frak{g}.v_\lambda]^{G_{v_\lambda}^\circ}_{- \lambda}
$$ car $T_{v_\lambda}$ agit trivialement sur le membre de droite.

Ainsi, on a $[V(\lambda)/ \frak{g}.v_\lambda]^{G_{v_\lambda}} =
[V(\lambda)/ \frak{g}.v_\lambda]^{\frak{g}_{v_\lambda}}_0 \oplus
[V(\lambda)/ \frak{g}.v_\lambda]^{\frak{g}_{v_\lambda}}_{- \lambda}
$

On en déduit alors la proposition. En effet, on a :

$$\frak{g}_{v_\lambda}=   \frak{u} \oplus \frak{t}_{v_\lambda}  \oplus  \displaystyle{ 
\bigoplus_{ \begin{array}{c} \alpha \in R_+, \,   \langle \lambda, \alpha ^ \vee  \rangle =0\end{array} }  \frak{g}_{- \alpha}
 } $$

Tout vecteur de poids $- \lambda$ est invariant par les algèbres de Lie $\frak{t}_{v_\lambda} $ et\\ $ \displaystyle{ 
\bigoplus_ { \begin{array}{c}     \alpha \in R_+, \,  \langle \lambda , \alpha ^\vee \rangle =0\end{array}    } \frak{g}_{- \alpha} } $, et tout vecteur de poids $0$ aussi s'il est invariant par $ \frak{u}$.
\hfill $\Box$

On obtient maintenant une condition nécessaire  pour que l'espace tangent soit non nul :

\begin{prop}\

Si $[V(\lambda)/ \frak{g}.v_\lambda]^{U}_0 \neq 0$, alors $\lambda$ s'écrit $\lambda = \alpha + \beta$ où $\alpha  \in S$ et $\beta \in R_+$.

Si $[V(\lambda)/ \frak{g}.v_\lambda]^{U}_{-\lambda} \neq 0$, alors $\lambda$ s'écrit $\lambda = \frac{\alpha + \beta}{2}$ où $\alpha  \in S$ et $\beta \in R_+$.

\end{prop}

\noindent {\bf Preuve.} Soit $v  \in V(\lambda)_0$ dont la classe dans  $V(\lambda)/ \frak{g}.v_\lambda$ est un $U$-invariant non nul. On exprime cela à l'aide de l'algèbre de Lie $\frak{u}$ de $U$ :

   $$\frak{u}v \subseteq \frak{g}v_\lambda \mbox{ et } v \notin \frak{g}v_\lambda.$$

Comme $v \notin V(\lambda)^U$, il existe une racine simple $\alpha$ telle que $ e_\alpha v \neq 0$.
Donc $ e_\alpha v$ est un $T$-vecteur propre de $\frak{g}.v_\lambda$.
Si  $ e_\alpha v$ était proportionnel à $v_\lambda$, alors $v \in V(\lambda)_{\lambda - \alpha}= \frak{g}_{- \alpha} . v_{\lambda}$ : une contradiction.
Donc il existe une racine positive $\beta$, telle que  $ e_\alpha v$ est proportionnel à  $ e_{-\beta} v_\lambda$. En considérant les poids, on obtient : $\alpha  = - \beta + \lambda $.

On vérifie de même la seconde implication.
\hfill $\Box$

\begin{prop}

Supposons l'espace tangent à $H_\lambda$ en $C_\lambda$ non nul, et $G$ semi-simple.
Alors l'image de $G$ dans $\operatorname{GL} (V(\lambda))$ est simple ou de type $A_1 \times A_1$. 
\end{prop}

\noindent {\bf Preuve.}

L'algèbre de Lie de $G$ est un produit d'algèbres de Lie simples : $ \frak{g} = \frak{g}_1 \times ... \times \frak{g}_r$, avec $r \geq 1$.
Celle de $U$ s'écrit : $ \frak{u} = \frak{u}_1 \times ... \times \frak{u}_r$, avec $ \frak{u}_j \subseteq \frak{g}_j$.

La donnée du poids dominant $\lambda$ de $\frak{g}$ revient à celle d'un poids dominant $\lambda_i$ de $\frak{g}_i$ pour tout $i$, et $$V(\lambda)=V(\lambda_1) \otimes_\Cbb ... \otimes_\Cbb V(\lambda_r).$$

On peut supposer que tous les $\lambda_i$ sont non nuls.

D'apres la proposition précédente, $\lambda$ est somme ou demi-somme de deux racines, on peut donc se limiter au cas où $r = 1 \mbox{ ou } 2$. Le cas $r=1$ correspond au cas où $G$ est simple~; supposons donc que $r=2$, et que l'espace tangent en $C_\lambda$ est non nul.

Selon la proposition 1.8, il existe un vecteur non nul $v$ appartenant à $ V(\lambda )_0 \mbox{ ou } V(\lambda )_{ -\lambda }$ qui est $U$-invariant modulo $\frak{g}.v_\lambda$.

Le vecteur $v$ n'est invariant ni par $\frak{u}_1$, ni par $\frak{u}_2$.
Donc il existe une racine simple $\alpha_1$ de $\frak{g}_1$ (resp. $\alpha_2$ de $\frak{g}_2$) telle que : $e_{\alpha_1}.v \not = 0$ (resp. $e_{\alpha_2}.v \not = 0$).
Comme $e_{\alpha_1}.v$ (resp. $e_{\alpha_2}.v$) est dans $\frak{g}.v_\lambda$, il existe une racine $\beta _1$ telle que $e_{\alpha_1}.v$ est proportionnel à $e_{-\beta_1}.v_\lambda$ (resp. une racine $\beta _2$ telle que $e_{\alpha_2}.v$ est proportionnel à $e_{-\beta_2}.v_\lambda$). 

Supposons que $v$ appartient à $ V(\lambda )_0$. En considérant les poids, on a $$\alpha_1 + \beta_1 = \alpha_2 + \beta_2 = \lambda,$$ donc $\beta_1$ est en fait une racine de $\frak{g}_2$, et $\beta_2$ une racine de $\frak{g}_1$. On a donc, avec des notations évidentes $$ (\alpha_1 , 0) +  ( 0 , \beta_1)= (\lambda _ 1 , \lambda _ 2)  \mbox{ et } ( 0 , \alpha_2) + (  \beta_2 , 0)=(\lambda _ 1 , \lambda _ 2).$$
Donc $$\lambda_1 = \alpha _ 1 \mbox{ et } \lambda_2 = \alpha _ 2.$$

Lorsque $v$ appartient à $V(\lambda )_{-\lambda}$ on en déduit de même $$\lambda_1 = \alpha _ 1 /2 \mbox{ et } \lambda_2  = \alpha _ 2 /2.$$

Dans les deux cas, les algèbres de Lie simples $\frak{g}_1$ et $\frak{g}_2$ admettent une racine simple qui est un poids dominant : elles sont donc de type $A_1$.
\hfill $\Box$

\subsubsection{Cas restant à étudier}\

Selon la proposition 1.9, pour que l'espace tangent en $C_\lambda$ soit non nul, il faut que $\lambda$ soit somme ou demi-somme d'une racine simple et d'une racine positive. On dresse ci-dessous la liste des cas où l'espace tangent peut être non nul, obtenue en calculant toutes les sommes et demi-sommes d'une racine simple et d'une racine positive, puis en ne gardant que celles qui sont des poids dominants. \

Les notations sont celles de \cite{Bou}.\

On donne $\lambda$  sous la forme d'une somme de poids fondamentaux, et sous la forme  de somme ou demi-somme d'une racine simple et d'une racine positive de toutes les fa{\c c}ons possibles.\\

\begin{enumerate}[(C1)]

\item $G$ est simple de type $A_1$,\\ $\lambda = 2 \omega _ 1 = \frac{1 }{2}( \alpha _ 1 + \alpha _ 1)$

\item $G$ est simple de type $A_1$,\\ $\lambda = 4 \omega _ 1 =  \alpha _ 1 + \alpha _ 1$

\item $G$ est simple de type $A_2$,\\  $\lambda = 3 \omega _ 1 =  \alpha _ 1 + ( \alpha _ 1 + \alpha _ 2 )$

\item $G$ est simple de type $A_2$,\\  $\lambda = 3 \omega _ 2 =  \alpha _ 2 + ( \alpha _ 1 + \alpha _ 2 )$

\item $G$ est simple de type $A_3$,\\  $\lambda =  \omega _ 2 = \frac{1 }{2}[ \alpha _ 2 + ( \alpha _ 1 + \alpha _ 2 + \alpha _ 3)] $

\item $G$ est simple de type $A_3$,\\ $\lambda = 2 \omega _ 2 =  \alpha _ 2 + ( \alpha _ 1 + \alpha _ 2 + \alpha _ 3)$

\item $G$ est simple de type $A_n$, $n \geq 2$,\\  $\lambda =  \omega _ 1 +  \omega _ n =  \alpha _ 1 + ( \alpha _ 2 +... + \alpha _ n ) =  \alpha _ n + ( \alpha _ 1 +  ... +  \alpha _ {n - 1})$

\item $G$ est simple de type $B_2$,\\  $\lambda =  \omega _ 2 = \frac{1}{2}[ \alpha _ 2 + (\alpha _ 1 +  \alpha _ 2 )] $

\item $G$ est simple de type $B_2$,\\  $\lambda =  2 \omega _ 2 = \alpha _ 2 + (\alpha _ 1 +  \alpha _ 2 ) $

\item $G$ est simple de type $B_3$,\\  $\lambda =  \omega _ 3 = \frac{1}{2}[ \alpha _ 3 + (\alpha _ 1 + 2 \alpha _ 2 + 2 \alpha _ 3)]$

\item $G$ est simple de type $B_3$,\\  $\lambda = 2 \omega _ 3 = \alpha _ 3 + (\alpha _ 1 + 2 \alpha _ 2 + 2 \alpha _ 3) $

\item $G$ est simple de type $B_n$, $n \geq 3$,\\ $\lambda =  \omega _ 2 = \alpha _ 2 + ( \alpha _ 1 + \alpha _ 2 + 2 ( \alpha _ 3 +...+  \alpha _ n)) $

\item $G$ est simple de type $B_n$, $n \geq 2$,\\  $\lambda =  \omega _ 1 = \alpha _ 1 + ( \alpha _ 2+...+  \alpha _ n)  = \alpha _ n + (\alpha _ 1 +...+\alpha _{n - 1} ) \\ = \frac{1}{2}[ \alpha _ 1 + (\alpha _ 1 +2(\alpha _ 2+...+  \alpha _ n ))] $ 

\item $G$ est simple de type $B_n$, $n \geq 2$,\\  $\lambda =  2 \omega _ 1 = \alpha _ 1 + (\alpha _ 1 +2(\alpha _ 2+...+  \alpha _ n )) $

\item $G$ est simple de type $C_n$, $n \geq 3$,\\  $\lambda =  \omega _ 1 = \frac{1}{2}[  \alpha _ 1 + (\alpha _ 1 +2(\alpha _ 2+...+ \alpha _{n - 1} ) + \alpha _ n) ] $

\item $G$ est simple de type $C_n$, $n \geq 3$,\\  $\lambda = 2 \omega _ 1 =  \alpha _ 1 + (\alpha _ 1 +2(\alpha _ 2+...+ \alpha _{n - 1} ) + \alpha _ n) $

\item $G$ est simple de type $C_n$, $n \geq 3$,\\  $\lambda =  \omega _ 2 = \alpha _ 1 +(2(\alpha _ 2+...+ \alpha _{n - 1} ) + \alpha _ n) \\  = \alpha _ 2 + (\alpha _ 1+ \alpha _ 2 +2(\alpha _ 3+...+ \alpha _{n - 1} ) + \alpha _ n)  $

\item $G$ est simple de type $D_4$,\\  $\lambda =  \omega _ 3 = \frac{1}{2}[  \alpha _ 3 +( \alpha _ 1+2 \alpha _2  + \alpha _ 3 + \alpha _ 4 )]  $ 

\item $G$ est simple de type $D_4$,\\  $\lambda =  \omega _ 4 = \frac{1}{2}[ \alpha _ 4 +( \alpha _ 1+2 \alpha _2  + \alpha _ 3 + \alpha _ 4)]   $

\item $G$ est simple de type $D_4$,\\  $\lambda = 2 \omega _ 3 =  \alpha _ 3 +( \alpha _ 1+2 \alpha _2  + \alpha _ 3 + \alpha _ 4) $ 

\item $G$ est simple de type $D_4$,\\  $\lambda = 2 \omega _ 4 = \alpha _ 4 +( \alpha _ 1+2 \alpha _2  + \alpha _ 3 + \alpha _ 4 ) $

\item $G$ est simple de type $D_n$, $n \geq 4$,\\  $\lambda =  \omega _ 1 = \frac{1}{2}[\alpha _ 1+( \alpha _ 1+ 2(\alpha _ 2 +...+  \alpha _ {n-2}) + \alpha _ {n-1} +   \alpha _ n)  ] $

\item $G$ est simple de type $D_n$, $n \geq 4$,\\  $\lambda = 2 \omega _ 1 = \alpha _ 1+( \alpha _ 1+ 2(\alpha _ 2 +...+  \alpha _ {n-2}) + \alpha _ {n-1} +   \alpha _ n)  $ 

\item $G$ est simple de type $D_n$, $n \geq 4$,\\  $\lambda =  \omega _ 2 =\alpha _ 2 +( \alpha _ 1 +\alpha _ 2 + 2(\alpha _ 3 +...+  \alpha _ {n-2}) +  \alpha _ {n-1}  +  \alpha _ n)  $ 
 
\item $G$ est simple de type $E_6$,\\  $\lambda = \omega _ 2= \alpha _ 2 + ( \alpha _ 1 +  \alpha _ 2 + 2 \alpha _ 3 + 3 \alpha _ 4+ 2 \alpha _ 5 +  \alpha _ 6 )  $

\item $G$ est simple de type $E_7$,\\  $\lambda =\omega _ 1=\alpha _ 1 + ( \alpha _ 1 + 2 \alpha _ 2 + 3  \alpha _ 3 + 4 \alpha _ 4+ 3 \alpha _ 5 + 2 \alpha _ 6 +  \alpha _ 7)   $

\item $G$ est simple de type $E_8$,\\  $\lambda =\omega _ 8=\alpha _ 8+( 2 \alpha _ 1 + 3 \alpha _ 2 + 4 \alpha _ 3 + 6 \alpha _ 4+ 5 \alpha _ 5 +  4 \alpha _ 6 + 3 \alpha _ 7 +  \alpha _ 8)  $

\item $G$ est simple de type $F_4$,\\  $\lambda =  \omega _ 1 =\alpha _ 1+(  \alpha _ 1 + 3 \alpha _ 2 + 4 \alpha _ 3 + 2 \alpha _ 4 ) $ 

\item $G$ est simple de type $F_4$,\\  $\lambda =  \omega _ 4 =  \alpha _ 3 + ( \alpha _ 1 + 2 \alpha _ 2 + 2 \alpha _ 3 + 2 \alpha _ 4) = \alpha _ 4 + ( \alpha _ 1 + 2 \alpha _ 2 + 3 \alpha _ 3 + \alpha _ 4 ) $ 

\item $G$ est simple de type $G_2$,\\  $\lambda =  \omega _ 1 = \alpha _ 1 + ( \alpha _ 1 +  \alpha _ 2)=  \frac{1}{2}[ \alpha _ 1 +  ( 3 \alpha _ 1 + 2 \alpha _ 2) ]  $ 

\item $G$ est simple de type $G_2$,\\  $\lambda = 2 \omega _ 1 =   \alpha _ 1 +  ( 3 \alpha _ 1 + 2 \alpha _ 2)   $

\item $G$ est simple de type $G_2$,\\  $\lambda =  \omega _ 2 =    \alpha _ 2 +  ( 3 \alpha _ 1 +  \alpha _ 2)   $

\item $G$ est semi-simple de type $A_1 \times A_1$,\\  $\lambda =  (  \omega_1 ,  \omega_1)=  \frac{1}{2}[( \alpha _ 1 , 0) + (0 , \alpha _ 1) ]     $

\item $G$ est semi-simple de type $A_1 \times A_1$,\\  $\lambda = ( 2 \omega_1 , 2 \omega_1)= ( \alpha _ 1 , 0) + (0 , \alpha _ 1)    $

\end{enumerate}

Tous les cas se traitent de la même manière : selon la proposition 1.8, il s'agit de calculer $$\dim [V(\lambda)/ \frak{g}.v_\lambda]^{U}_0 \mbox{ et }
\dim [V(\lambda)/ \frak{g}.v_\lambda]^{U}_{- \lambda}.$$
Des connaissances élémentaires sur les modules irréductibles (pour lesquelles on renvoie à \cite{Se}) permettent de mener à bien le calcul. A titre d'exemples, on traite quelques cas représentatifs dans les sections 1.3.2 à 1.3.4.

\subsubsection{ Cas de la représentation adjointe d'un groupe simple}

Il s'agit des cas (C1), (C7), (C9), (C12), (C16), (C25), (C26), (C27), (C28) et (C32).

On a $V(\lambda) \cong  \frak{g}$ et $\lambda$ est la plus grande racine.\\

1) Le cas (C1) où $G$ est de type $A_1$ se traite à part : $ \frak{g}$ s'identifie à l'algèbre de Lie des matrices $2 \times 2$ de trace nulle. 
On pose :
 
\[x := \begin{pmatrix}
         0 & 1 \\
         0 & 0
\end{pmatrix} \qquad  
h := \begin{pmatrix}
         1 & 0 \\
         0 & -1
\end{pmatrix} \qquad
y := \begin{pmatrix}
         0 & 0 \\
         1 & 0
\end{pmatrix}\]

La matrice $x$ est un vecteur de plus grand poids : l'espace tangent est isomorphe à $$[\frak{g}/ \frak{g}.x]^U_{- \lambda} = [\frak{g}/ ( \Cbb x \oplus \Cbb h ) ]^U_{- \lambda} \cong  \Cbb y.$$

Il est donc de dimension $1$.\\

2) Dans les autres cas, $\lambda$ n'est pas demi-somme d'une racine simple et d'une racine positive. L'espace tangent est donc isomorphe à 
$[\frak{g}/ \frak{g}.v_\lambda]^U_0$.

Posons $$ E  := \{ \gamma \in R \mbox{ $\vert$ } \exists \delta \in R,\gamma + \delta = \lambda  \}.  $$

On a alors, en notant $h_\lambda$ l'élément de $[ \frak{g}_\lambda ,  \frak{g}_{-\lambda } ] $ tel que $\lambda (h_\lambda)=2$ : 
$$ \frak{g} v_\lambda = [ \frak{g},  \frak{g}_\lambda]= \frak{g}_\lambda \oplus \Cbb h_\lambda \oplus \bigoplus_ {\gamma \in E} \frak{g}_\gamma . $$

Le sous-espace de poids nul de $\frak{g}/ \frak{g}.v_\lambda$ est isomorphe à $ \frak{t} / \Cbb h_\lambda$, et l'espace tangent est isomorphe à $$ \{ t \in \frak{t} \mbox{ $\vert$ } \forall \alpha \in S,  [ \frak{g}_\alpha , \frak{t} ] \subseteq \frak{g} v_\lambda  \} / \Cbb h_\lambda ,$$donc à $$ \{ t \in \frak{t} \mbox{ $\vert$ } \forall \alpha \in S \setminus E ,  \alpha(t) = 0  \} / \Cbb h_\lambda .$$
Il est donc de dimension $ \dim \frak{t} - \operatorname{card} (S \setminus E) - 1$. Or $ \dim \frak{t} = \operatorname{card} S $. La dimension de l'espace tangent est donc $ \operatorname{card} ( S \cap E ) -1$.

Comme $ \operatorname{card} ( S \cap E ) $ est le nombre de fa{\c c}ons d'écrire $\lambda$ comme somme d'une racine simple et d'une racine positive, on conclut à l'aide de la liste 1.3.4 que l'espace tangent est de dimension $1$ dans le cas où $G$ est de type $A_n$, et $0$ dans les autres cas.

\subsubsection{ Cas (C23) }
Ici, $\frak{g}$ s'identifie à l'algèbre de Lie $\frak{so}(2n)$ des matrices de taille $2n \times 2n$ antisymétriques par rapport à la seconde diagonale.
On a une action naturelle de $\frak{so}(2n)$ sur $\Cbb ^ {2n}$, dont la base canonique est notée $(e_1, ... , e_n, e_{-n}, ... , e_{-1})$.
Le module simple $V(2 \omega_1)$ peut être vu comme un quotient du carré symétrique de $\Cbb ^ {2n}$ : $$V(2 \omega_1)= S^2 \Cbb ^ {2n}/ \langle e_1e_{-1}+ ... + e_n e_{-n} \rangle.$$
Il s'agit de calculer la dimension de $ [V(2 \omega_1)/\frak{g}v_{2 \omega_1}]^U_0$.\\

Notons $\pi  : V(2 \omega_1) \longrightarrow V(2 \omega_1)/\frak{g}v_{2 \omega_1} $ la surjection canonique.
On remarque que $\pi$ induit des isomorphismes :$$ \pi( V(2 \omega_1)_0) \cong V(2 \omega_1)_0, \mbox{ } \pi( V(2 \omega_1)_{\alpha_2}) \cong V(2 \omega_1)_{\alpha_2}, \mbox{ ... , }  \pi( V(2 \omega_1)_{\alpha_n}) \cong V(2 \omega_1)_{\alpha_n}.$$

Par contre, $\pi( V(2 \omega_1)_{\alpha_1}) = 0$.
On a donc, en notant $\langle A \rangle$ la sous-algèbre de $\frak{g}$ engendrée par une partie $A$ de $\frak{g}$ :

 $$ [\pi (V(2 \omega_1))]^U_0 \cong V(2 \omega_1)_0^{ \langle  e_ { \alpha _2},...,e_{\alpha_n} \rangle  }=V(2 \omega_1)_0^{ \langle  e_ { \alpha _2},...,e_{\alpha_{n-1}} \rangle  }$$
car $e_ { \alpha _{n-1}} \mbox{ et } e_{\alpha_{n}} $ stabilisent les mêmes éléments dans $V(2 \omega_1)_0$.

Or on sait que $V(2 \omega_1)_0$ est de dimension $n-1$, et $V(2 \omega_1)_{ \alpha _{1}}, ... ,V(2 \omega_1)_{ \alpha _{n-1}}$ sont de dimension $1$.\\

Comme  $V(2 \omega_1)_0^{ \langle  e_ { \alpha _1},...,e_{\alpha_{n-1}} \rangle  }=V(2 \omega_1)_0^U=\{0\}$, on en déduit que l'espace vectoriel $V(2 \omega_1)_0^{ \langle  e_ { \alpha _2},...,e_{\alpha_{n-1}} \rangle  }$ est de dimension $1$.

Donc l'espace tangent est de dimension $1$.

\subsubsection{ Cas (C31) }

Les poids de $V(2 \omega_1)$ sont :
\begin{itemize}
 \item $0$ avec multiplicité $3$
 \item $\alpha _ 1 $ et ses conjugués sous l'action du groupe de Weyl, chacun avec multiplicité $2$
 \item $\alpha _ 2 $ et ses conjugués, chacun avec multiplicité $1$
 \item $2 \omega_1 = 4 \alpha _ 1 + 2 \alpha _ 2 $ et ses conjugués, chacun avec multiplicité $1$.
 \end{itemize}

Comme $V(2 \omega_1)_0^U=\{0\}$, l'application linéaire
$$\begin{array}{rlc} \phi  : V(2 \omega_1)_0& \rightarrow& V(2 \omega_1)_{ \alpha _{1}} \oplus V(2 \omega_1)_{ \alpha _{2}}\\ v &\mapsto &(  e_ { \alpha _1 } v ,   e_ { \alpha _2} v)
\end{array}$$ est injective, donc bijective.

Notons $\pi  : V(2 \omega_1) \longrightarrow V(2 \omega_1)/\frak{g} . v_{2 \omega_1} $ la surjection canonique.
On remarque que $\pi$ induit des isomorphismes $$ \pi( V(2 \omega_1)_0) \cong V(2 \omega_1)_0 \mbox{ et } \pi( V(2 \omega_1)_{\alpha_2}) \cong V(2 \omega_1)_{\alpha_2}.$$ Par contre, $\dim \pi( V(2 \omega_1)_{\alpha_1}) = 1$.

On en déduit que l'application quotient $$ \bar { \phi }  : \pi( V(2 \omega_1)_0) \rightarrow \pi( V(2 \omega_1)_{ \alpha _{1}}) \oplus \pi (V(2 \omega_1)_{ \alpha _{2}})$$ a pour noyau $ \pi (V(2 \omega_1)_0)^U$, de dimension $1$.

Donc l'espace tangent est de dimension $1$.

\subsection{Conclusion}

Dans ce paragraphe, on ramène la démonstration du théorème 1.1, à des résultats obtenus au \S 2.4 (de fa{\c c}on indépendante). On commence par énoncer sans démonstration une conséquence immédiate du lemme de Nakayama :

\begin{lem}
Soit $A = \bigoplus _ { d \in \Nbb} A _ d $ un anneau gradué par $\Nbb$. On suppose l'anneau $A_0$ local d'idéal maximal $ \Mcal _ 0 $. Notons $ \Mcal$ l'idéal maximal $ \Mcal _ 0 \oplus \bigoplus _ {d =1}^\infty A_d$ de $A$.
Soit $M = \bigoplus _ {d \in \Zbb } M _ d$ un $A$-module gradué de type fini.
Si $ \Mcal . M = M$, alors $M=0$.
\end{lem}

Le schéma $H_\lambda$ est presque déterminé par le corollaire suivant du corollaire 1.5 et de la proposition 1.6 :

\begin{cor}
Le schéma $H_\lambda$ est soit une droite affine, soit un point épaissi $\Spec \Cbb [ t ]/(t^N)$, pour un entier $N$.
\end{cor}

\noindent {\bf Preuve.} On garde les notations du corollaire 1.5.
L'idéal $ \Mcal $ de $A$ correspondant au point $C_\lambda$ de $H_\lambda$ est l'unique idéal maximal de $A$ fixé par $\Gbb_m$, donc on a $$ \Mcal= \Mcal _ 0 \oplus \bigoplus _ {d =1}^\infty A _ d$$ où $ \Mcal _ 0$ est l'idéal maximal de $A_0$.
Selon la proposition 1.6, l'espace vectoriel $\Mcal / \Mcal^2$ est de dimension inférieure ou égale à 1. Il existe donc un élément homogène $f$ de $\Mcal$ tel que $\Mcal = A f + \Mcal^2 $. Selon le lemme précédent, $\Mcal = A f$.
Or $A$ s'écrit, comme espace vectoriel $A = \Cbb \oplus \Mcal$, donc on a $A=\Cbb \oplus \Cbb f  \oplus \Cbb f ^2 \oplus ... = \Cbb [ f ]$.
Ainsi, l'algèbre $A$ est monogène, et comme son spectre est connexe (proposition 1.3), on en déduit le corollaire.
\hfill $\Box$\\

Le théorème 1.1 est donc démontré dans les cas où l'espace tangent est nul.
Dans les cas où il est de dimension 1, il reste à exclure les points épaissis. On va voir (\S 2.4) qu'à chaque fois que l'espace tangent à $H_\lambda$ en $C_\lambda$ est de dimension 1, il existe d'autres points fermés que $C_\lambda$ dans $H_\lambda$, et celui-ci est donc une droite affine.

\section{Algèbres de Jordan simples et familles universelles}

On va maintenant relier la classification \textit{(H)} du théorème 1.1 à deux classifications déjà connues :
\begin{itemize}
 \item celle notée \textit{(J)} des algèbres de Jordan simples (théorème 2.1)
 \item et une classification notée \textit{(A)} de variétés projectives à deux orbites (théorème 2.2).
 \end{itemize}
Pour cela, on va associer de fa{\c c}on naturelle \begin{itemize}
 \item aux objets de \textit{(J)} des objets de \textit{(A)} dans \S 2.3 (en considérant le cône des éléments de rang $1$ des algèbres de Jordan simples)
 \item et aux objets de \textit{(A)} des objets de \textit{(H)} dans \S 2.4 
 \end{itemize}

On constatera alors que lors de ces deux opérations, tous les cas sont atteints.
Ainsi, pour chacun des cas du théorème 1.1, on obtiendra à l'aide d'une algèbre de Jordan simple une déformation non triviale du cône $C_\lambda$ dans $V(\lambda)$, qui sera en fait la famille universelle au dessus de $H_\lambda$ (\S 2.5).
\subsection{Classification des algèbres de Jordan simples }\

Pour plus de détails concernant les algèbres de Jordan, on renvoie à \cite{Jac} ou \cite{FK}.
On appellera \textit{algèbre de Jordan} (complexe) une $\Cbb$-algèbre  $(A, *)$ commutative unitaire de dimension finie (non nécessairemant associative) telle que $$\forall a,b \in A, ~ a^2 * ( a * b) = a * ( a^2 * b).$$
On peut montrer que $A$ est associative relativement aux puissances (c'est-à-dire : toutes ses sous-algèbres monogènes sont associatives).

On peut alors définir le polynôme minimal d'un élément $a$ de $A$ : c'est le générateur unitaire de l'idéal $\{ P \in \Cbb [X] \mbox{ tels que } P(a)=0 \}$ de $\Cbb [X]$.
Un élément de $A$ est dit \textit{régulier} si le degré de son polynôme minimal est maximal. Les éléments réguliers forment un ouvert dense de $A$.

Il existe alors (\cite{FK}, prop II.2.1) des fonctions polynômiales $p_1, ... , p_r$ sur $A$ telles que si $a$ est un élément régulier de $A$, son polynôme minimal est $X^r + p_1(a) X^{r-1}+ ... + p_r (a)$. Chaque $p_i$ est alors homogène de degré $i$.
On définit la trace et le déterminant sur $A$ par :$$\tr :=-p_1 \mbox{ et } \det :=(-1)^r p_r.$$

Dans la suite, on s'intéresse aux algèbres de Jordan \textit{simples}, que l'on définit maintenant.
Si $a$ est un élément de $A$, on note $L(a) :A \rightarrow A$, $b \mapsto a * b$ la multiplication par $a$. On note $\Tr$ la trace d'un endomorphisme $A \rightarrow A$.
Une algèbre de Jordan $A$ est \textit{semi-simple} si la forme bilinéaire
$$\begin{array}{cll}  A \times A& \rightarrow & \Cbb\\ (a,b) &\mapsto & \Tr L(a * b)
\end{array}$$ est non dégénérée. Elle est dite \textit{simple} si de plus elle n'admet pas d'idéaux non triviaux. (Dans ce cas, les formes linéaires $a \rightarrow \operatorname{Tr} L (a)$ et $a \rightarrow tr (a)$ sont en fait proportionnelles.)

Les matrices hermitiennes sur les complexifiées $$R= \Cbb \otimes_\Rbb  \Rbb \mbox{, } \Cbb \otimes_\Rbb \Cbb \mbox{, } \Cbb \otimes_\Rbb \Hbb \mbox{, } \Cbb \otimes_\Rbb \Obb$$ des algèbres de Hurwitz donnent des exemples d'algèbres de Jordan simples : on note $H_n(R)$ l'espace vectoriel des matrices de taille $n \times n$ qui sont égales à la transposée de leur conjuguée. On munit $H_n(R)$ d'une structure d'algèbre en posant :$$M_1 * M_2 =  \frac{1 }{2} ( M_1 M_2 +M_2 M_1).$$

La classification des algèbres de Jordan simples est connue (\cite{Jac}, théorème 8 p 203) :
\begin{teo}[P.Jordan, J.von Neumann, E.Wigner]
Toute algèbre de Jordan simple est isomorphe à l'une des suivantes :

(J1) $\Cbb \oplus W$, où $W$ est un espace vectoriel de dimension finie muni d'une forme bilinéaire non dégénérée $\langle.,.\rangle$. La loi de l'algèbre est donnée par la formule : $ (t_1,w_1)*(t_2,w_2)= (t_1 t_2 + \langle w_1,w_2 \rangle , t_1 w_2+t_2 w_1).$

(J2) $H_n(\Cbb \otimes_\Rbb \Rbb)$ , $n \geq 3$

(J3) $H_n(\Cbb \otimes_\Rbb \Cbb)$ , $n \geq 3$

(J4) $H_n(\Cbb \otimes_\Rbb \Hbb)$ , $n \geq 3$

(J5) $H_3(\Cbb \otimes_\Rbb \Obb)$.
\end{teo}

On définit enfin le \textit{groupe de structure} d'une algèbre de Jordan simple $A$ : c'est le groupe des automorphismes (d'espace vectoriel) de $A$ qui conservent à un scalaire près le déterminant :$$ \operatorname{Str}(A) :=\{ g \in GL(A) \mbox{ $\vert$ } \exists u \in \Cbb^*, \forall a \in A, \det (ga)=u \det(a) \}.$$
Il contient le groupe $\operatorname{Aut}(A)$ des automorphismes d'algèbre de $A$. Les deux groupes $\operatorname{Aut}(A)$ et $ \operatorname{Str}(A)$ sont des groupes algébriques réductifs (mais non connexes). En fait, $\operatorname{Aut}(A)^\circ$ est semi-simple, et $ \operatorname{Str}(A)^\circ$ est de centre les homothéties.

\subsection{Classification de variétés à deux orbites }\

Les espaces homogènes sous l'action d'un groupe réductif admettant une complétion équivariante par un diviseur homogène ont été classifiés par D. Akhiezer (voir \cite{Ak1} ou \cite{HuS}~; la classification est retrouvée dans \cite{Bri} par des méthodes algébriques). Dans le cas où le diviseur est ample, on a :

\begin{teo}[D.Akhiezer] Soit $Z$ une variété projective lisse.
Soit $D$ un diviseur ample de $Z$, et $\Omega$ son complémentaire. On suppose qu'il existe une action régulière d'un groupe algébrique affine connexe $\Gamma$ sur $Z$ sous laquelle $\Omega$ et $D$ sont les orbites de $Z$.

Soit $G$ l'image de $\Gamma$ dans $\operatorname{Aut}(Z)$, et $H$ le stabilisateur d'un point de $\Omega$. Alors, à revêtement fini de $G$ près, on est dans un des cas suivants :

(A1) $G=SL(n+1), n \geq 1$ , $ H=GL(n)$ et $Z = \Pbb ^ n  \times (\Pbb ^ n)^*$.

(A2) $G=SO(n), n \geq 3$ , $ H=SO(n-1)$ et $Z = Q(n-1)$.

(A3) $G=SO(n), n \geq 3$ , $ H=O(n-1)$ et $Z = \Pbb ^ {n-1}$.

(A4) $G=Sp(2n), n \geq 2$ , $ H=Sp(2) \times Sp(2n-2)$ et $Z$ est la grassmanienne des 2-plans de $\Cbb ^ {2n}$.

(A5) $G=F_4$ , $H= Spin(9)$ et $Z= E_6 / P$ en notant $P$ le sous-groupe parabolique maximal de $E_6$ dont les racines simples sont $ \alpha_2, ... , \alpha_n $ (notations de \cite{Bou}).

(A6) $G=G_2$ , $ H=SL(3)$ et $Z=Q(6)$.

(A7) $G=G_2$ , $ H=N_G(SL(3))$ et $Z=\Pbb ^ 6 $.

(A8) $G=Spin(7)$ , $ H=G_2$ et $ Z = Q(7)$.

(A9) $G=SO(7)$ , $ H=G_2$ et $Z=\Pbb ^ 7 $.\\
On a noté $Q(n)$ une quadrique projective lisse de dimension $n$.\\
Les actions de $SL(n+1)$, $SO(n)$ et $Sp(2n)$ sont les actions naturelles.
Les actions de $G_2$, $Spin(7)$ et $SO(7)$ sont déduites des plongements $G_2 \hookrightarrow SO(7)$ via la représentation de dimension 7 de $G_2$, et $Spin(7) \hookrightarrow SO(8)$ via la représentation spinorielle de $Spin(7)$. 
\end{teo}

\subsection{ Un lien entre les classifications (J) et (A) }\
 
Dans ce paragraphe, on rappelle une définition du cône des éléménts \textit{de rang $1$} d'une algèbre de Jordan simple. La variété projective formée des droites de ce cône nous donne alors une variété à 2 orbites de la classification \textit{(A)}~; l'orbite ouverte de cette variété correspond aux éléments de trace non nulle.

On obtient au passage une bijection entre les classes d'isomorphisme des \textit{algèbres de Jordan simples complexes} et les classes d'isomorphisme des \textit{espaces symétriques de rang 1 complexes} (en effet, ceux-ci sont les quotients $G/H$, où $G$ et $H$ sont les groupes donnés dans les cas \textit{(A1)} à \textit{(A5)}).

La correspondance analogue dans le cas réel (entre les \textit{algèbres de Jordan simples réelles} et les \textit{espaces symétriques de rang 1 réels compacts}) est établie dans \cite{Hirz}~; dans ce cas, tout élément de rang 1 est de trace \textit{non nulle}, c'est pourquoi l'espace symétrique des droites d'éléments de rang 1 (et de trace non nulle) est \textit{compact}.\\

Soit $A$ une algèbre de Jordan simple.

Notons $\Gamma$ la composante neutre du groupe des automorphismes de $A$, et $\Gamma^\prime$ celle du groupe de structure de $A$ :$$\Gamma :=\mbox{Aut}(A)^\circ \subseteq \Gamma^\prime  := \mbox{Str}(A)^\circ .$$
Notons $\Cbb.1$ la droite engendrée par l'élément unité de $A$, et $V$ le sous-espace vectoriel de $A$ formé des éléments de trace nulle. Alors $A$ est somme directe des deux sous-espaces vectoriels :
\begin{equation}\label{d}
 A=\Cbb.1 \oplus V,    
\end{equation}
et cette décomposition est stable par $\Gamma$.
On vérifie (grâce à la classification) que $V$ est un $\Gamma$-module simple.

Notons $D$ la $\Gamma$-orbite fermée dans $\Pbb(V)$ et $ \widetilde{D}$ le cône des vecteurs primitifs de $V$, c'est-à-dire le cône affine sur $D$.

On vérifie également que comme $\Gamma^\prime $-module rationnel, $A$ est simple~; notons $Z$ la $\Gamma^\prime$-orbite fermée dans $\Pbb(A)$ et $\widetilde{Z}$ son cône des vecteurs primitifs de $A$.
Les éléments (non nuls) de $\widetilde{Z}$ sont appelés \textit{les éléments de rang $1$} de l'algèbre de Jordan $A$. Les éléments (non nuls) de $\widetilde{D}$ sont les éléments de rang $1$ et de trace nulle.

On remarque que $D$ est un diviseur ample de $Z$ :$$D = Z \cap \Pbb (V) \subseteq Z \subseteq \Pbb(\Cbb.1 \oplus V)$$et on vérifie enfin que $\Gamma$ agit transitivement sur $D$ et sur $Z \setminus D$.

Ainsi, à toute algèbre de Jordan simple on fait correspondre un élément de la classification (A) en prenant comme groupe $G$ le groupe $\Gamma$~; on peut aussi prendre comme groupe $G$ un sous-groupe fermé de $\Gamma$ pourvu qu'il agisse transitivement sur $D$ et sur $Z \setminus D$.

Précisément :

\begin{itemize}
\item à partir de \textit{(J1)}, on obtient le cas \textit{(A2)} quand $G$ est le groupe $\Gamma = SO(W)$, mais aussi, en considérant des sous-groupes stricts de $\Gamma$, le cas \textit{(A6)} quand $W$ est de dimension $7$ et $G=G_2$ et le cas \textit{(A8)} quand  $W$ est de dimension $8$ et $G=\operatorname{Spin}(7)$.
 \item à partir de \textit{(J2)}, on obtient le cas \textit{(A3)} quand $G$ est le groupe $\Gamma= SO(n)$, mais aussi, en considérant des sous-groupes stricts de $\Gamma$, les cas \textit{(A7)} quand $n=7$ et \textit{(A9)} quand $n=8$.
 \item à partir de \textit{(J3)}, on obtient le cas \textit{(A1)} avec $G= \Gamma = \operatorname{PGL}(n)$.
 \item à partir de \textit{(J4)}, on obtient le cas \textit{(A4)} avec $G= \Gamma = \operatorname{Sp}(2n)$.
 \item à partir de \textit{(J5)}, on obtient le cas \textit{(A5)} avec $G= \Gamma = F_4$.
 \end{itemize}

On voit donc que tous les cas du théorème 2.2 peuvent être obtenus à partir des algèbres de Jordan simples.\\

\subsection{ Un lien entre les classifications (A) et (H) }\

Supposons que le groupe réductif connexe $G$ agit sur une variété projective lisse $Z$, et que ses orbites sont un diviseur ample $D$ et son complémentaire $\Omega$ (de sorte que l'on est dans la situation du théorème 2.2).

Alors $D$ est en fait très ample, et si l'on plonge $Z$ dans $\Pbb( \Gamma ( Z , \Ocal  _Z (D))^*)$ en associant à tout $z \in Z$ l'hyperplan des sections globales qui s'annulent en $z$, le cône affine au-dessus de $Z$ dans $\Gamma ( Z , \Ocal  _Z (D))^*$ est normal (car selon le \S2.3, c'est le cône des vecteurs primitifs d'un $G$-module simple).
Ce cône affine est donc le spectre de l'algèbre graduée $$R  := \bigoplus _ {d \in \Nbb} \Gamma ( Z , \Ocal  _Z (dD)),$$ où l'on note $ \Ocal  _Z (dD)$ le faisceau inversible sur $Z$ associé au diviseur $dD$, et $ \Gamma ( Z , \Ocal  _Z (dD))= :R_d$ l'espace de ses sections globales.
On note $\sigma _ D \in R _ 1$ la section canonique de $\Ocal  _Z (D)$.
L'algèbre $R$ est naturellement munie d'une structure de $G$-algèbre rationnelle~; on munit $\widetilde{Z}$ de l'action de $G$ correspondante (qui induit celle de $G$ sur $Z= \operatorname{Proj}R$).

\begin{prop}

Il existe un poids dominant $\lambda$ tel que $R_1$ se décompose comme $G$-module sous la forme $$R_1 = \Cbb \sigma _ D \oplus V ( \lambda )^ *.$$ Notons $f$ l'immersion fermée correspondant au morphisme surjectif d'algèbres $\operatorname{Sym} (\Cbb \sigma _ D \oplus V(\lambda)^*) \longrightarrow R$ et $\pi$ le morphisme donné par la fonction régulière $\sigma _ D$ : on a un diagramme commutatif de morphismes équivariants

$$\xymatrix{ \widetilde Z  \ar[rr]^{f~~~} \ar[rd]_{\pi} && \Abb ^ 1 \times V ( \lambda ) \ar[ld] ^ {pr_1} \\ & \Abb ^ 1 } $$
où $\Abb ^ 1$ est muni de l'action triviale de $G$.
De plus $\pi  : \widetilde Z \longrightarrow \Abb ^ 1$ est une famille de fonction de Hilbert $h _ \lambda$ et la fibre de $\pi$ en $0 \in \Abb ^ 1 $ est le cône des vecteurs primitifs de $V(\lambda)$.
Les autres fibres $\pi ^ {-1} ( t ) $, $t \neq 0$ sont isomorphes à l'orbite ouverte $\Omega$.
\end{prop}

\noindent {\bf Preuve.}

Comme $D$ est complet et homogène sous l'action de $G$, il est isomorphe à un quotient $G/P$, où $P$ est un sous-groupe parabolique de $G$ contenant $B$.
On a donc un plongement $i : G / P   \hookrightarrow Z$.
L'image réciproque de $\Ocal  _Z (D)$ par $i$ est un faisceau inversible ample sur $G/P$ : elle est donc isomorphe au faisceau $\Lcal _ \lambda$ pour un certain poids dominant $\lambda$.

On a une suite exacte de $\Ocal  _Z$-modules :$$0 \longrightarrow \Ocal  _Z (-D) \overset{\sigma _ D}{\longrightarrow} \Ocal  _Z  \longrightarrow i _ * \Ocal _ { G / P} \longrightarrow 0 .$$
On la tensorise par $\Ocal  _Z (dD)$ :$$0 \longrightarrow \Ocal  _Z ((d-1)D) \overset{\sigma _ D}{\longrightarrow}  \Ocal  _Z (dD) \longrightarrow i _ * \Lcal _ { d \lambda} \longrightarrow 0 .$$ 
On a donc une suite exacte de $G$-modules de sections globales :

\begin{equation}\label{S}
0 \longrightarrow R _ {d - 1 } \overset{\sigma _ D}{\longrightarrow}  R _ d \overset{ f _ d}{\longrightarrow} V( d \lambda ^ *).
\end{equation}
Lorsque $d=1$, la suite (5) est $$0 \longrightarrow \Cbb \overset{\sigma _ D}{\longrightarrow}  R _ 1 \overset{ f _ 1}{\longrightarrow} V( \lambda ^ *).$$ Comme $Z$ se plonge dans $\Pbb ( R_1^*)$, l'espace vectoriel $R_1$ n'est pas de dimension $1$. De plus $V( \lambda) ^ *$ est un $G$-module simple, donc le morphisme $f_1$ est surjectif, et on en déduit le premier point de la proposition.

Montrons que le morphisme $f_d$ est surjectif pour tout entier $d$. Comme $R_1$ contient un $B$-vecteur propre de poids $\lambda^*$, $R_d$ contient un $B$-vecteur propre de poids $d \lambda^*$, et il contient donc un $G$-module simple isomorphe à  $V( d \lambda ^ *)$. Or en considérant la suite exacte (5) pour tout $d^ \prime < d$, on remarque que $R_{d-1}$ est un sous-$G$-module de $\bigoplus _ {d=0} ^ {d-1} V( d \lambda ^ *)$, d'où l'assertion.

(On peut retrouver la surjectivité des $f_d$ par un argument cohomologique. En effet, selon le \S 2.3, la variété $Z$ est une variété de drapeaux pour l'action d'un groupe $G^\prime$ réductif connexe, que l'on peut supposer simplement connexe quitte à le remplacer par un revêtement fini. L'algèbre des fonctions régulières sur $G^\prime$ est alors factorielle, et le groupe de Picard de $G^\prime$ est nul. Selon \cite{KKV}, prop 3.2 (i), tout faisceau inversible sur $Z$ est donc linéarisable. On peut donc appliquer le théorème de Borel-Weil-Bott (\cite{Ak2} p 113 ) au faisceau inversible $\Ocal  _Z ((d-1)D)$. Comme celui-ci est ample, on obtient $H^1(Z,\Ocal  _Z ((d-1)D))=0$, d'où le résultat.)

La fibre de $\pi$ au dessus de $0$ est le sous-cône de $\widetilde{Z}$ d'algèbre affine graduée  $R/{\sigma _ D}R $. D'après ce qui précède, on a un isomorphisme de $G$-modules $$R/{\sigma _ D}R \simeq \bigoplus_{d \in \Nbb}  V( d \lambda ^ *).$$

La fibre au dessus de $0$ est donc le cône des vecteurs primitifs de $V(  \lambda)$, selon la proposition 1.3(a).

La fibre au dessus de $t \neq 0$ est la section de $Z$ par l'hyperplan affine $\{ \sigma _ D = t \}$, donc est isomorphe à l'ouvert $\Omega$.
Enfin le morphisme $\pi$ est plat car $R$ est un $\Cbb [\sigma _ D ]$-module sans torsion : le morphisme $\pi$ est donc bien une famille de fonction de Hilbert $h_\lambda$.
\hfill $\Box$

\begin{rem} La déformation $\pi$ ainsi obtenue est toujours non triviale, car les fibres $\{ \pi ^ {-1} ( t ), t \not = 0 \}$ sont homogènes pour l'action de $G$, contrairement à $ \pi ^ {-1} ( 0 )$.
\end{rem}

On associe ainsi à chaque objet de \textit{(A)} un schéma de Hilbert invariant de la classification \textit{(H)}, et on constate que l'on obtient ainsi toute cette classification :

\begin{itemize}
 \item Le cas \textit{(H2)} provient du cas \textit{(A1)}, avec $n \geq 2$.
 \item Le cas \textit{(H1)} (resp. \textit{(H4)}, \textit{(H6)}, \textit{(H9)}) provient des cas \textit{(A2)} et \textit{(A3)}, avec $n=3$ (resp. $n \text{ impair supérieur à } 5$, $n \text{ pair supérieur à } 6$, $n=4$).
 \item Le cas \textit{(H3)} provient des cas \textit{(A8)} et \textit{(A9)}.
 \item Le cas \textit{(H5)} provient du cas \textit{(A4)}, avec $n \geq 3$.
 \item Le cas \textit{(H7)} provient du cas \textit{(A5)}.
 \item Le cas \textit{(H8)} provient des cas \textit{(A6)} et \textit{(A7)}.
\end{itemize}

En particulier, on en déduit dans chacun des cas l'existence d'un point de $H_\lambda$ distinct de $C_\lambda$, comme annoncé au \S 1.4.

\subsection{ Construction des familles universelles}\

Pla{\c c}ons-nous dans l'un des cas du théorème 1.1 : le schéma de Hilbert invariant $H_\lambda$ est isomorphe à la droite affine. Il résulte des \S 2.3 et 2.4 qu'on obtient une déformation de $C_\lambda$ à partir d'une (unique) algèbre de Jordan simple $A$, de la fa{\c c}on suivante. On note $\widetilde Z$ le cône des éléments de rang $1$ de $A$, et $i$ l'inclusion $\widetilde Z \subseteq A$.
Le morphisme donné par la restriction de la trace de $A$ à $ \widetilde Z$ est noté $\pi_\lambda$.
Le diagramme analogue à celui de la proposition 2.3 est alors :
$$\xymatrix{ \widetilde Z  \ar[rr]^{i~~~~~~~~~} \ar[rd]_{\pi_\lambda} && A \simeq  \Abb ^ 1 \times V ( \lambda ) \ar[ld] ^ {pr_1} \\ & \Abb ^ 1 } $$

\begin{prop}
Si $A$ n'est pas de type \textit{(J1)}, la famille $\pi_\lambda$ est la famille universelle au dessus de $H_\lambda \simeq \Abb^1$.\\
Sinon on a $A = \Cbb \oplus W$~; la famille universelle est alors, avec des notations évidentes $$\xymatrix{ \{(t,w) \in \Abb^1 \times W ~ \vert ~ t=\langle w,w\rangle \}  \ar[rr]^{~~~~~~~~~~~~~i} \ar[rd]_{t} &&  \Abb ^ 1 \times W \ar[ld] ^ {pr_1} \\ & \Abb ^ 1 } $$
\end{prop}

\noindent {\bf Preuve.}
Les familles de la proposition sont les images inverses par un (unique) morphisme $f  : \Abb ^1 \rightarrow H _ \lambda$ de la famille universelle (car ce sont bien des familles de fonction de Hilbert $h_\lambda$).

On va montrer que $f$ est injectif~; comme $H_\lambda$ est isomorphe à $ \Abb ^1$, on en conclura que $f$ est un isomorphisme, et la proposition sera démontrée.

L'injectivité de $f$ signifie que les fibres des familles de sous-schémas considérées sont deux à deux distinctes.

Cela est clair dans le cas où $A$ est de type \textit{(J1)} : la fibre de la famille au dessus de $t \in \Abb ^1$ est la sous-variété $\{w \in W ~ \vert ~ t= \langle w,w\rangle \ \}$.

Dans le cas où $A$ n'est pas de type \textit{(J1)}, on vérifie que son élément unité n'est pas somme de deux éléments de rang $1$.\\
Le module simple $V( \lambda )$ est l'espace $V$ de la décomposition $(2)$. La fibre de $\pi _ \lambda$ au dessus de $t \in \Abb ^1$ est $$\{a-t.1 ~ \vert ~ a \in \widetilde Z \mbox{ et } \operatorname{tr}(a)=t  \} \subseteq V.$$
Supposons que les fibres de $\pi _ \lambda$ au dessus de $t$ et $t ^ \prime$ soient égales : on peut alors écrire $ a-t.1 = a^ \prime - t^ \prime .1 $, donc $(t-t^ \prime).1=a^ \prime-a$, donc $t=t^ \prime$.

Le morphisme $f$ est donc bien injectif.
\hfill $\Box$

Le second cas de la proposition correspond aux cas \textit{(H1)},\textit{(H4)},\textit{(H6)},\textit{(H9)} où $G=\operatorname{SO}(W)$ et $V( \lambda )=W$, ainsi qu'au cas \textit{(H8)} où $G=G_2$ et $V( \lambda )=W$ est de dimension $7$, et au cas \textit{(H3)} où $G=\operatorname{Spin}(7)$ et $V( \lambda )=W$ est de dimension $8$.

\section{Rigidité des cônes de vecteurs primitifs}

On sait (\cite{H}, ex 9.8 p 267) que les déformations infinitésimales de $C_\lambda$ sont classifiées par un $\Cbb$-espace vectoriel noté $T^1(C_\lambda)$. Comme $C_\lambda$ est une $G$-variété, l'espace vectoriel $T^1(C_\lambda)$ est un $G$-module rationnel (\cite{Rim})~; on le note dans la suite $T^1_\lambda$.

On voit facilement que l'espace des éléments $G$-invariants de $T^1_\lambda$ est en fait l'espace tangent en $C_\lambda$ au schéma de Hilbert invariant (proposition \ref{t3.4})~; il est donc déterminé par le théorème 1.1. Dans cette partie, on détermine complètement le $G$-module $T^1_\lambda$ :
\begin{teo}

L'espace $T^1_\lambda$ des déformations infinitésimales de $C_\lambda$ est nul, sauf dans les cas suivants :
\begin{enumerate}[(R1)]
 \item Si l'on est dans les cas (H2) à (H9) du théorème 1.1, alors $T^1_\lambda=V(0)$.
 \item Si $G = \operatorname{SL}(2)$ et $m \geq 2$ est un entier, alors $T^1_{m}=V(m-2) \oplus V(m-4)$ (on indexe les poids de $\operatorname{SL}(2)$ par les entiers).
 \item Si le groupe s'écrit $G=\operatorname{SL}(2) \times H$ et le module simple $V(\lambda)=V_{\operatorname{SL}(2)}(m) \otimes W, $ où $V_{\operatorname{SL}(2)}(m)$ est le $\operatorname{SL}(2)$-module simple de plus grand poids $m$, et $(H,W)= ( \operatorname{SL}(V),V) $, $ ( \operatorname{SL}(V),V^*)$ ou $ (\operatorname{Sp}(V),V)$ (dans ce dernier cas, $V$ est un espace vectoriel de dimension paire supérieure ou égale à $2$), alors $T^1_\lambda=V_{\operatorname{SL}(2)}(m-2) \otimes W$.
 \end{enumerate}

\noindent et dans les cas obtenus à partir d'un cas parmi les précédents par factorisation.
\end{teo}

\begin{rems}

\begin{enumerate}[(1)]

 \item
On retrouve ainsi des faits déjà connus :
\begin{itemize}
 \item Le cône affine sur le plongement de Segre de la variété $\Pbb^m \times \Pbb^n$ (dans $\Pbb ^ {(m+1)(n+1)-1}$) est rigide quand $m+n \geq 3$ (voir par exemple \cite{KL} thm 2.2.8).
 \item Pinkham a déterminé dans \cite{Pink} l'espace $T^1_m$ quand $C_m$ est le cône affine sur la courbe rationnelle normale de degré $m$ dans $\Pbb ^m$~; en particulier, il a montré que $\operatorname{dim} (T^1_m)=2m-4$ (pour $m \geq 4$). Il a aussi montré que si $m \geq 5$, la déformation verselle est irréductible, de dimension $m-1$, et lisse hors de l'origine. Si $m=4$, elle a deux composantes de dimensions $3$ et $1$ qui se rencontrent transversalement à l'origine.

Le cône de vecteurs primitifs $C_4$ est exceptionnel, car c'est le seul dont l'espace des déforma\-tions infinitésimales admet à la fois une partie $G$-invariante et une partie non invariante : $T^1_4 = V(0) \oplus V(2)$. La direction $G$-invariante de l'espace $T^1_4$ correspond à la composante de dimension $1$ de la déformation verselle.

Ainsi, on constate que dans tous les cas où le schéma de Hilbert invariant $H_\lambda$ n'est pas réduit à un point, il donne une composante irréductible de la déformation verselle de $C_\lambda$.
 \item Svanes a montré dans \cite{Sv1} et \cite{Sv2} que les cônes affines sur les plus petits plongements des variétés de drapeaux de $\operatorname{SL}(n)$ (qui correspondent au cas où $G$ est simple de type $A_n$ et $\lambda$ est une somme de poids fondamentaux) sont rigides, à l'exception des cas \textit{(H2)} et, pour $n=3$, \textit{(H6)} du théorème 1.1.
\end{itemize}

 \item
Les couples $(H,W)$ du troisième cas du théorème peuvent être décrits géométriquement : ce sont ceux où le groupe $H$ agit transitivement sur les droites du module $W$ (cela résulte par exemple de \cite{Ak2} thm2 p75).

\end{enumerate}
\end{rems}

Enfin on détermine les déformations verselles des c\^ones $C_\lambda$. Si on est
dans le  cas \textit{(R1)} du théorème 3.1, la déformation verselle de
$C_\lambda$ est donnée par le schéma de Hilbert invariant $H_\lambda$. Dans le
cas \textit{(R2)}, elle a été déterminée dans \cite{Pink} par des équations
explicites.

Plaçons-nous maintenant dans le cas \textit{(R3)} du théorème 3.1. On va
décrire la déformation verselle de
$C_\lambda$  à l'aide d'équations
analogues à celles de \cite{Pink}. La dimension de
l'espace vectoriel $W$ est $n+1$, pour un entier non nul $n$.
Le c\^one $C_\lambda$ est le c\^one affine dans $\Abb^{(m+1)(n+1)}$ au dessus
de l'image de $\Pbb^1 \times \Pbb^n$ par le plongement
 $$\begin{array}{cll} \Pbb^1 \times \Pbb^n & \longrightarrow&
 \Pbb^{(m+1)(n+1)-1}\\ ([a:b],[c_0,...,c_n]) &\longmapsto &
 [a^jb^{m-j}c_i]_{_{i,j}}
\end{array}$$
On note ici ce c\^one $C_{mn}$. Son idéal homogène est engendré par les
mineurs $2 \times 2$ de la matrice  $$\displaystyle{\left(\begin{array}{cccc|cccc|c|cccc} x_{00} & x_{01} &  \cdot \cdot \cdot &
      x_{0~m-1} &  x_{10} & x_{11} & \cdot \cdot \cdot  & x_{1~m-1} &   \cdot \cdot \cdot &  x_{n0} & x_{n1} & \cdot \cdot \cdot & x_{n~m-1} \\ x_{01} & x_{02} &  \cdot \cdot \cdot &
      x_{0~m} &  x_{11} & x_{12} &  \cdot \cdot \cdot & x_{1~m} &  \cdot \cdot \cdot &  x_{n1} & x_{n2} &  \cdot \cdot \cdot & x_{n~m}  \end{array}\right) }.$$

\begin{prop}
La déformation verselle de $C_{mn}$ est la déformation $\frak{V}$ au dessus du spectre de
l'anneau $\Cbb[[\mathbf{t}]]$ des séries formelles en les $t_{i,j}$ (où $0
\leq i  \leq n$ et $ 1 \leq  j \leq m-1$) définie par les
mineurs $2 \times 2$ de la matrice
\begin{equation}\label{verse} \left(\begin{array}{cccc|c|cccc} x_{00} &
      \cdot \cdot \cdot &  x_{0~m-2} &
      x_{0~m-1} &   \cdot \cdot \cdot &  x_{n0} & \cdot \cdot \cdot & x_{n~m-2} & x_{n~m-1} \\ x_{01}-t_{01} &  \cdot \cdot \cdot  &  x_{0~m-1} - t_{0~m-1} &
      x_{0~m} &  \cdot \cdot \cdot &  x_{n1} -t_{n1} & \cdot \cdot \cdot & x_{n~m-1} - t_{n~m-1} & x_{n~m}  \end{array}\right)            \end{equation}

\end{prop}

On démontre le théorème 3.1 dans les parties 3.1 à 3.3, et la proposition 3.3
dans la partie 3.4. On rappelle d'abord quelques faits connus.\\

Notons $\Tcal_{C_\lambda}$ et $\Tcal_{V(\lambda)}$ les faisceaux tangents respectifs de $C_\lambda$ et $V(\lambda)$, et $\Ncal_{C_\lambda}$ le faisceau normal de $C_\lambda$ dans $V(\lambda)$.
On a $\Tcal_{V(\lambda)}= \Ocal_{V(\lambda)} \otimes V(\lambda)$. Comme $C_\lambda$ est normal, $\Tcal_{C_\lambda}$ et $\Ncal_{C_\lambda}$ sont des faisceaux réflexifs.

Les déformations infinitésimales de $C_\lambda$ se plongent en fait toutes dans $V(\lambda)$, et l'on a une suite exacte (\cite{H}, ex 9.8 p 267)
$$0 \longrightarrow H^0(C_\lambda, \Tcal_{C_\lambda}) \longrightarrow  H^0(C_\lambda, \Tcal_{V(\lambda)} |_{C_\lambda}) \longrightarrow H^0(C_\lambda, \Ncal_{C_\lambda}) \longrightarrow T^1_\lambda \longrightarrow 0 ,$$ c'est-à-dire
\begin{equation}\label{ty6}
0 \longrightarrow H^0(C_\lambda, \Tcal_{C_\lambda}) \longrightarrow  H^0(C_\lambda, \Ocal_{C_\lambda}) \otimes_\Cbb V(\lambda)  \longrightarrow H^0(C_\lambda, \Ncal_{C_\lambda}) \longrightarrow T^1_\lambda \longrightarrow 0 .   
\end{equation}

On peut supposer $C_\lambda$ de dimension supérieure ou égale à $2$~; on note $E_\lambda  := C_\lambda \setminus \{0\}$ le cône épointé. La suite exacte ci-dessus s'identifie alors à la suivante (avec des notations analogues) \begin{equation}\label{ty7} 0 \longrightarrow H^0(E_\lambda, \Tcal_{E_\lambda}) \longrightarrow  H^0(E_\lambda, \Ocal_{E_\lambda}) \otimes_\Cbb V(\lambda)  \longrightarrow H^0(E_\lambda, \Ncal_{E_\lambda}) \longrightarrow T^1_\lambda \longrightarrow 0 .\end{equation}

Or comme $E_\lambda$ est lisse, on a la suite exacte courte $$0 \longrightarrow \Tcal_{E_\lambda} \longrightarrow \Ocal_{E_\lambda} \otimes_\Cbb V(\lambda) \longrightarrow \Ncal_{E_\lambda} \longrightarrow 0 .$$

On en déduit la proposition suivante, due à Schlessinger (\cite{Schle}) :

\begin{prop}\label{t3.3}
On a une suite exacte : $$0 \longrightarrow T^1_\lambda \longrightarrow H^1(E_\lambda, \Tcal_{E_\lambda}) \longrightarrow  H^1(E_\lambda,  \Ocal_{E_\lambda}) \otimes_\Cbb V(\lambda).$$
\end{prop}

\subsection{Préliminaires}

On commence par déterminer la partie invariante de l'espace $T^1_\lambda$ :

\begin{prop}\label{t3.4} On a un isomorphisme canonique $$(T^1_\lambda)^G \cong T_{C_\lambda} H_\lambda.$$ Ainsi, l'espace  $(T^1_\lambda)^G$ est nul, sauf dans les cas (H1) à (H9) du théorème 1.1, où il est de dimension 1.

\end{prop}

\noindent {\bf Preuve.}

En prenant les $G$-invariants de (\ref{ty6}), on obtient la suite exacte de $G$-modules de \cite{AB}, prop 1.13 : $$ 0 \longrightarrow  H^0(C_\lambda, \Tcal_{C_\lambda})^G  \longrightarrow ( H^0(C_\lambda, \Ocal_{C_\lambda})\otimes V(\lambda))^G  \longrightarrow T_{C_\lambda} H_\lambda \longrightarrow (T^1_\lambda)^G  \longrightarrow 0, $$ qui s'écrit dans notre cas (\cite{AB}, prop 1.15 (iii)) :  $$ 0 \longrightarrow [\frak{g}.v_\lambda]^{G_{v_\lambda}}  \longrightarrow V(\lambda)^{G_{v_\lambda}} \longrightarrow T_{C_\lambda} H_\lambda \longrightarrow (T^1_\lambda)^G  \longrightarrow 0. $$
Comme ses deux premiers termes sont de dimension $1$ : $$ [\frak{g}.v_\lambda]^{G_{v_\lambda}} = V(\lambda)^{G_{v_\lambda}} = \Cbb v_\lambda,$$ on en déduit le résultat.
\hfill $\Box$

Ainsi, on a montré que la partie $G$-invariante était bien celle annoncée dans le théorème 3.1. On va maintenant déterminer les autres composantes isotypiques de $T^1_\lambda$.\\

On note $X_\lambda$ la variété de drapeaux $G/P_\lambda$, et  $\pi  : E_\lambda \longrightarrow X_\lambda$ la surjection naturelle. On remarque que $E_\lambda$ est l'espace total du faisceau $\Lcal_\lambda$ privé de la section nulle, donc $\pi$ est un morphisme affine lisse.

\begin{prop}\label{t3.5} On a une suite exacte de faisceaux $G$-linéarisés sur $X_\lambda$ : $$0  \longrightarrow \bigoplus_{d \in \Zbb} \Lcal_{d \lambda} \longrightarrow \pi_* \Tcal_{E_\lambda} \longrightarrow \bigoplus_{d \in \Zbb} \Lcal_{d \lambda} \otimes \Tcal_{X_\lambda} \longrightarrow 0, $$ donc une suite exacte de $G$-modules :
\begin{equation}\label{ty8}
\bigoplus_{d \in \Zbb} H^1(X_\lambda , \Lcal_{d \lambda} )  \rightarrow  H^1(E_\lambda,\Tcal_{E_\lambda})  \rightarrow  \bigoplus_{d \in \Zbb} H^1(X_\lambda , \Lcal_{d \lambda} \otimes \Tcal_{X_\lambda}) \rightarrow  \bigoplus_{d \in \Zbb} H^2(X_\lambda , \Lcal_{d \lambda} ).   
\end{equation}
\end{prop}

\noindent {\bf Preuve.}

Comme le morphisme $\pi$ est lisse, on a la suite exacte courte $$ 0 \longrightarrow \Tcal_\pi \longrightarrow \Tcal_{E_\lambda} \longrightarrow \pi^* \Tcal_{X_\lambda} \longrightarrow 0.$$ On remarque que le faisceau  $\Tcal_\pi$ tangent à $\pi$ est isomorphe à $\Ocal_{E_\lambda}$.
Puis, comme $\pi$ est affine, on a $$ 0 \longrightarrow \pi_* \Tcal_\pi \longrightarrow \pi_* \Tcal_{E_\lambda} \longrightarrow \pi_* \pi^* \Tcal_{X_\lambda} \longrightarrow 0,$$ d'où la suite exacte de faisceaux annoncée, car $\pi_* \Tcal_\pi \simeq \bigoplus_{d \in \Zbb} \Lcal_{d \lambda}$, et, selon la formule de projection, $\pi_* \pi^* \Tcal_{X_\lambda} \simeq  \bigoplus_{d \in \Zbb} \Lcal_{d \lambda} \otimes \Tcal_{X_\lambda}$.

La suite exacte de $G$-modules donnée en découle aussi, car, comme $\pi$ est affine,\\ $H^1(E_\lambda,\Tcal_{E_\lambda})  \cong H^1(X_\lambda, \pi_* \Tcal_{E_\lambda} )$.
\hfill $\Box$

La suite exacte (\ref{ty8}) nous permettra de démontrer le théorème 3.1 dans certains cas, à l'aide du théorème de Borel-Weil-Bott (\cite{Ak2}, thm p113).

Notons $Q_\lambda$ l'unique sous-groupe parabolique de $G$ conjugué à $P_\lambda$ et contenant le sous-groupe de Borel $B^-$ opposé à $B$~; on voit naturellement $Q_\lambda$ comme un point de $X_\lambda= G/P_\lambda$.
Regardons, pour appliquer le théorème de Borel-Weil-Bott, quels sont les poids de l'action de $T$ sur les fibres $\Lcal_{d \lambda}|_{ \{ Q_\lambda \} }$ et $\Tcal_{X_\lambda} |_{ \{ Q_\lambda \} } \cong \frak{g}/ { {\frak{q}}_\lambda}$ : le tore $T$ agit sur le premier espace avec le poids $d \lambda ^ *$~; les poids du second sont les racines positives de $G$ qui ne sont pas des racines de $Q_\lambda$.

Si $W$ est un $Q_\lambda$-module rationnel, on note $\Vcal(W)$ le $\Ocal_{X_\lambda}$-module $G$-linéarisé dont la fibre en $Q_\lambda$ est $W$. On rappelle (voir \cite{Gross} ou \cite{jant}) que l'espace des sections globales de $\Vcal(W)$ est le $G$-module induit par le $Q_\lambda$-module $W$ : $$ H^0(X_\lambda ,\Vcal(W)) = \operatorname{Ind}_ {Q_\lambda} ^ G ( W )$$ et les groupes de cohomologie de $\Vcal(W)$ donnent les foncteurs dérivés à droite du foncteur $ \operatorname{Ind}_ {Q_\lambda} ^ G$ : $$ H^j(X_\lambda ,\Vcal(W)) = R^j \operatorname{Ind}_ {Q_\lambda } ^ G ( W ).$$ 
Dans toute la suite, on note simplement $ \operatorname{Ind}$ le foncteur $\operatorname{Ind}_ {Q_\lambda} ^ G$.\\
Si $\mu$ est un caractère de $Q_\lambda$, on définit le $Q_\lambda$-module $$W[\mu] :=\Cbb_\mu \otimes_\Cbb W,$$ où $\Cbb_\mu$ est la droite où $Q_\lambda$ opère avec le poids $\mu$.

On note enfin $\star$ l'action tordue du groupe de Weyl sur $\Lambda$ : si $w$ est un élément du groupe de Weyl, et $\mu$ un poids, on pose $w \star \mu  := w(\mu + \rho) - \rho$, où $\rho$ est la demi-somme des racines positives.

\begin{prop}\label{t3.6} Soit $d \in \Zbb$. Si l'espace $H^1(X_\lambda , \Lcal_{d \lambda} )$ est non nul, alors $V(\lambda)$ est en fait un $\operatorname{SL}(2)$-module, et l'action de $G$ se factorise sous la forme $G \longrightarrow \operatorname{SL}(2) \longrightarrow \operatorname{GL}(V(\lambda)).$
\end{prop}

\noindent {\bf Preuve.} Si $d \geq 0$, on sait que tous les groupes de cohomologie de $\Lcal_{d \lambda}$ sont nuls sauf en degré $0$.

Supposons $d<0$. Selon le théorème de Borel-Weil-Bott, $H^1(X_\lambda , \Lcal_{d \lambda} )$ est nul ou irréductible, et on a $H^1(X_\lambda , \Lcal_{d \lambda} ) \simeq V(\mu)$ si et seulement si il existe une racine simple $\alpha$ et un  poids dominant $\mu$ tels que $$ s _ \alpha \star \mu = d \lambda ^ *,$$ en notant $ s _ \alpha$ la réflexion simple associée à $ \alpha$.
On en déduit $$ \mu - d \lambda ^ * = ( 1 + \langle \alpha ^ \vee , \mu \rangle ) \alpha .$$
La racine simple $\alpha$ est donc un poids dominant : c'est une racine simple de $G$ isolée dans son diagramme de Dynkin.  Notons $\omega _ \alpha = \alpha /2$ le poids fondamental associé à $\alpha$. Le poids $-d \lambda ^ *$ est proportionnel à $\omega _ \alpha$, et $\lambda = \lambda ^ *$ aussi, d'où le résultat.
\hfill $\Box$

\begin{prop}\label{t3.7} Lorsque $d \geq 0$, on a $H^1(X_\lambda , \Lcal_{d \lambda} \otimes \Tcal_{X_\lambda})=0$.
\end{prop}

\noindent {\bf Preuve.} C'est une conséquence du fait suivant (\textit{cf.} \cite{Bro}, thm2.2) : notons $T^* X_\lambda$ le fibré cotangent de $X_\lambda = G/P_\lambda$, et $p : T^* X_\lambda \longrightarrow X_\lambda$ la projection canonique.
On a, pour tout $i\geq 1$, $$H^i(T^*X_\lambda , p^* \Lcal_{d \lambda} )=0.$$
Selon la formule de projection, $$f_*f^* \Lcal_{d \lambda} = \Lcal_{d \lambda} \otimes f_* \Ocal_{T^* X_\lambda} = \Lcal_{d \lambda} \otimes \operatorname{Sym}(\Tcal_{X_\lambda}),$$ où $\operatorname{Sym}(\Tcal_{X_\lambda})= \bigoplus _{k \in \Nbb} \operatorname{S}^k \Tcal_{X_\lambda}$ est l'algèbre symétrique du $\Ocal_{X_\lambda}$-module $\Tcal_{X_\lambda}$.
Ainsi, $$\bigoplus _{k \in \Nbb} H^i(X_\lambda , \Lcal_{d \lambda} \otimes \operatorname{S}^k \Tcal_{X_\lambda})=0 ,$$ et la proposition en découle en prenant $i=k=1$.
\hfill $\Box$

\subsection{Cas où le groupe $G$ est simple}

Dans ce paragraphe, on établit le théorème 3.1 dans le cas où $G$ est un groupe simple.

Le \S3.2.1 concerne le cas où $G=\operatorname{SL}(2)$. L'espace des déformations infinitésimales a alors été déterminé dans \cite{Pink}~; on donne cependant une preuve simple de ce résultat, qui a l'avantage de fournir la structure de $\operatorname{SL}(2)$-module de $T^1_\lambda$. 

Le cas des autres groupes simples est traité dans le \S3.2.2.

\subsubsection{Cas où $G=\operatorname{SL}(2)$}

On note $H$ le sous-groupe des matrices unipotentes triangulaires inférieures $\displaystyle{\left(\begin{array}{cc} 1 & 0\\ x & 1 \end{array}\right) }$.
On identifie le groupe des poids de $\operatorname{SL}(2)$ à $\Zbb$, de sorte que les poids dominants sont les éléments de $\Nbb$. Le poids dominant $\lambda$ est donc un entier, que l'on note ici $m$. On suppose $m \geq 1$~; on a donc $Q_\lambda=B^-$.

Pour déterminer $T^1 _ \lambda$, on utilise la suite exacte (\ref{ty7}).
Comme $\pi$ est affine, on a la suite exacte
$$0 \rightarrow H^0(X_\lambda, \pi _ * \Tcal_{E_\lambda}) \xrightarrow{f}  H^0(X_\lambda, \pi _ * \Ocal_{E_\lambda}) \otimes V(\lambda)  \xrightarrow{g} H^0(X_\lambda, \pi _ * \Ncal_{E_\lambda}) \rightarrow T^1_\lambda \rightarrow 0 .$$

On remarque que les fibres respectives en $Q_\lambda$ des faisceaux $G$-linéarisés $\pi _ * \Tcal_{E_\lambda}$, $\pi _ * \Ocal_{E_\lambda} \otimes V(\lambda)$  et $\pi _ * \Ncal_{E_\lambda}$ sont les $Q_\lambda$-modules gradués $$                                                                                                                                                     \bigoplus _{d \in \Zbb}  \frak{g}v_{ - \lambda}[md]                                                                      \mbox{ , }                                                                             \bigoplus _{d \in \Zbb}  V(m)  [md]                                                                                             \mbox{ et }                                                                            \bigoplus _{d \in \Zbb}  V(m)/  \frak{g}v_{ - \lambda}[md]     ,        $$ et que les morphismes $f$ et $g$ sont les images par le foncteur $\operatorname{Ind}$ des morphismes de modules gradués qui forment en chaque degré $d$ une suite exacte courte : $$0 \longrightarrow \frak{g}v_{ - \lambda}[md] \longrightarrow  V(m)  [md]  \longrightarrow  V(m)/  \frak{g}v_{ - \lambda}[md] \longrightarrow  0 .$$

Considérons la suite exacte longue associée à cette dernière
\begin{equation}\label{ty9} 0 \rightarrow \operatorname{Ind}(\frak{g}v_{ - \lambda}[md]) \xrightarrow{f_d} \operatorname{Ind}(  V(m)  [md]) \xrightarrow{g_d}  \operatorname{Ind}( V(m)/  \frak{g}v_{ - \lambda}[md])   \rightarrow R^1 \operatorname{Ind}(\frak{g}v_{ - \lambda}[md]) \rightarrow  ... \end{equation}

On a $$T^1_\lambda \cong  \bigoplus _{d \in \Zbb} \operatorname{coker} g _ d.$$

Pour connaître $\operatorname{coker} g _ d$, on remarque qu'on a des isomorphismes de $Q_\lambda$-modules :
$$\frak{g}v_{ - \lambda} \simeq V(1)[-m+1] \mbox{ et } V(m)/\frak{g}v_{ - \lambda}               \simeq V(m-2)[2],  $$ donc
$$\frak{g}v_{ - \lambda}[md] \simeq V(1)[m(d-1)+1]$$ et $$ V(m)/\frak{g}v_{ - \lambda}[md]               \simeq V(m-2)[md+2]. $$ La suite exacte longue (\ref{ty9}) devient alors                                                                                                                                                                                                                                                                                  $$0 \rightarrow  V(1) \otimes V(m(d-1)+1) \xrightarrow{f_d} V(m) \otimes V(md)  \xrightarrow{g_d}  V(m-2) \otimes V(md+2) $$ $$ \rightarrow R^1 \operatorname{Ind} ( V(1)[(m(d-1)+1])\rightarrow ...$$

On peut alors conclure :

\noindent Supposons $m \geq 3$ :
\begin{itemize}
 \item si $d<0$, on a  $V(m-2) \otimes V(md+2)=0$, donc $\operatorname{coker} g _ d =0$.
 \item si $d \geq 1$, montrons que $R^1 \operatorname{Ind} ( V(1)[m(d-1)+1])=0$ : considérons la suite exacte de $Q_\lambda$-modules suivante $$0 \longrightarrow \Cbb [m(d-1)] \longrightarrow V(1) [m(d-1)+1]  \longrightarrow \Cbb [m(d-1)+2] \longrightarrow 0 .$$ La suite exacte longue associée (relativement au foncteur $\operatorname{Ind}$) s'écrit $$... \longrightarrow 0 \longrightarrow R^1 \operatorname{Ind} (V(1) [m(d-1)+1])  \longrightarrow 0 \longrightarrow ...$$ Ainsi le morphisme $g_d$ est surjectif :  $\operatorname{coker}  g _ d =0$.

 \item si enfin $d=0$, on remarque que $\operatorname{coker} g _ 0 = V(m-2) \oplus V(m-4)$.
 \end{itemize}

\noindent Si $m=2$, le conoyau de $g_d$ est nul si $d \neq -1$, et on a $\operatorname{coker} g _ {-1}=V(0)$.
Donc on a toujours $T^1_\lambda = V(m-2) \oplus V(m-4)$.
\subsubsection{Autres groupes simples}

On suppose maintenant que $G$ est un groupe simple de type autre que $A_1$, et il s'agit de montrer que les seules déformations infinitésimales de $C_\lambda$ sont celles qui proviennent du théorème 1.1.

Selon la proposition  \ref{t3.4}, il ne reste qu'à montrer que le groupe $G$ agit trivialement sur l'espace $T^1_\lambda$  (\textit{ie} tous les éléments de l'espace sont $G$-invariants).

En vertu de la proposition \ref{t3.6}, la suite exacte (\ref{ty8}) donne ici $$  H^1(E_\lambda,\Tcal_{E_\lambda})  \hookrightarrow  \bigoplus_{d \in \Zbb} H^1(X_\lambda , \Lcal_{d \lambda} \otimes \Tcal_{X_\lambda}).$$ Il ne reste donc plus qu'à montrer que pour tout $d$ et pour tout poids dominant $\mu$ non nul, la composante isotypique $H^1(X_\lambda , \Lcal_{d \lambda} \otimes \Tcal_{X_\lambda})_{(\mu)}$ est nulle.

Pour cela on va utiliser le lemme et la proposition suivants.
Les notations sont celles de \cite{Bou}.

\begin{lem}\label{t3.8} Soit $R$ un système de racines irréductible muni d'une base $S$.
Soit $\alpha$ un élément de $S$.
Alors il existe une racine positive longue $\gamma$ telle que $\langle \gamma^\vee, \alpha \rangle = -1$, sauf dans les cas suivants :

\begin{itemize}
 \item si $R$ est de type $A_1$.
 \item si $R$ est de type $B_2$ et $\alpha = \alpha _1$ est la racine simple longue.
 \item si $R$ est de type $C_n$, avec $n\geq 3$ et $\alpha = \alpha _n$ est la racine simple longue.
 \end{itemize}
Lorsqu'une telle racine $\gamma$ existe, elle n'est pas unique, sauf dans les cas suivants :
\begin{itemize}
 \item si $R$ est de type $A_2$.
 \item si $R$ est de type $B_2$ et $\alpha = \alpha _2$ est la racine simple courte : seule $\gamma =\alpha _1$ convient. 
 \item si $R$ est de type $C_n$, avec $n\geq 3$ et $\alpha = \alpha _i$, $i=1...n-1$ est une racine simple courte : seule $\gamma = 2 \alpha _{i+1}+...+ 2 \alpha _{n-1}+ \alpha _n$ convient .
 \item si $R$ est de type $G_2$ : pour $\alpha = \alpha _1$, seule $\gamma =\alpha _2$ convient~; pour $\alpha = \alpha _2$, seule $\gamma =3 \alpha_1+ \alpha _2$ convient.
 \end{itemize}
\end{lem}

\noindent {\bf Preuve.}
Traitons d'abord le cas où toutes les racines de $R$ sont de même longueur. Si $R$ est de type $A_1$, il n'y a rien à prouver~; on suppose donc que $R$ est de type $A_n $ $ (n\geq 2)$, $D_n$ $ (n\geq 4)$, $E_6$, $E_7$ ou $E_8$. 
L'existence de $\gamma$ est claire : toute racine simple reliée à $\alpha$ dans le diagramme de Dynkin de $R$ convient. Montrons que si $R$ n'est pas de type $A_2$, $\gamma$ n'est pas unique.
Cela est clair si $\alpha$ est reliée à plusieurs racines simples dans le diagramme de Dynkin.
Sinon, notons $\beta$ la seule racine simple reliée à $\alpha$. Le sous-système de racines $R^\prime$ de $R$ de base $S^\prime  := S \setminus \{ \alpha \}$ est un système de racines irréductible, de rang supérieur ou égal à $2$. On sait qu'il admet plusieurs racines $\gamma$ telles que $\beta$ a pour coefficient $1$ dans l'écriture de $\gamma$ dans la base $S^\prime$ (voir par exemple \cite{Ak2}, prop 1 p 126). Toutes ces racines $\gamma$ conviennent clairement.

Supposons maintenant que $R$ est de type $B_n$, avec $n \geq 3$.
Si $\alpha$ est une racine simple longue, on est ramené au premier cas, car les racines longues de $R$ forment un système de racines de type $A_3$ si $n=3$, et $D_n$ sinon.
Sinon, $\alpha = \alpha _ n$, et la racine $\gamma = \alpha _ i +...+ \alpha _ {n-1}$ convient pour tout $i=1...n-1$.

Supposons que $R$ est de type $F_4$.
Si $\alpha$ est une racine simple longue, on est ramené au premier cas, car les racines longues de $R$ forment un système de racines de type $D_4$.
Si $\alpha = \alpha _ 3= \epsilon_4,$ alors $\gamma=\epsilon_i - \epsilon_4$ convient pour tout $i=1,2,3.$ 
Si $\alpha = \alpha _ 4= (\epsilon_1-\epsilon_2-\epsilon_3-\epsilon_4)/2,$ alors $\gamma=\epsilon_i + \epsilon_j$ convient pour tout $2 \leq i < j \leq 4$.

Enfin, on vérifie aisément les assertions du lemme concernant $B_2$, $C_n$ ($ n \geq 3$) et $G_2$.
\hfill $\Box$

\begin{prop}\label{t3.9} Soit $R$ un système de racines irréductible muni d'une base.
On suppose que $R$ n'est pas de type $A_1$.

Soient $\alpha$ une racine simple, $\beta$ une racine positive, et $N \geq 2$ un entier, tels que $N \alpha + \beta$ est un poids dominant.
Alors $N=2$, et on est dans l'un des cas suivants :
\begin{itemize}
\item si $R$ est de type $A_2$, on a $2 \alpha_1 + \alpha_2= 3 \omega_1$ et $2 \alpha_2 + \alpha_1= 3 \omega_2$.
 \item si $R$ est de type $B_2$, on a $2 \alpha_2 + \alpha_1= 2 \omega_2$.
 \item si $R$ est de type $C_n$, on a $2 \alpha_1 + (2 \alpha_2+ 2 \alpha_3+...+2 \alpha _{n-1}+ \alpha _n)= 2 \omega_1$.
 \item si $R$ est de type $G_2$, on a $2 \alpha_1 + \alpha_2= 3 \omega_1$.
 \end{itemize}

\end{prop}

\noindent {\bf Preuve.}
Traitons d'abord le cas où il existe plusieurs racines positives longues $\gamma$ telles que $\langle \gamma ^\vee  ,  \alpha \rangle = -1$.
Soit $\gamma$ une telle racine, que l'on suppose distincte de $\beta$.
On a $\langle \gamma ^\vee , N \alpha + \beta \rangle = -N + \langle \gamma ^\vee , \beta \rangle \geq 0.$
D'où $N \leq \langle \gamma ^\vee , \beta \rangle \leq 1$, car $\gamma$ est une racine longue distincte de $\beta$.

Pour conclure, on étudie un à un les cas du lemme précédent où il n'existe pas de racine $\gamma$, ainsi que ceux où il existe une unique racine $\gamma$ (selon le premier point de la démonstration, on peut alors supposer $\beta = \gamma$).
\hfill $\Box$

\begin{prop}\label{t3.10} Le groupe $G$ agit trivialement sur $H^1(X_\lambda , \Lcal_{d \lambda} \otimes \Tcal_{X_\lambda})$.
\end{prop}

\noindent {\bf Preuve.} Selon la proposition \ref{t3.7}, on peut supposer $d<0$.

Afin d'appliquer le théorème de Borel-Weil-Bott, on considère une suite de Jordan-Hölder du $Q_\lambda$-module $\frak{g}/ { {\frak{q}}_\lambda}$, c'est-à-dire une suite décroissante $$\frak{g}/ { {\frak{q}}_\lambda} = W_0 \supset W_1 \supset ... \supset W_r = 0$$ de $Q_\lambda $-modules telle que les quotients $W_i/W_{i+1}$ sont des modules simples.

Soit $\mu$ un poids dominant non nul.
On veut montrer que la composante isotypique $$H^1(X_\lambda , \Lcal_{d \lambda} \otimes \Tcal_{X_\lambda})_{(\mu)}= R^1 \operatorname{Ind} (\frak{g}/ { {\frak{q}}_\lambda} [d \lambda^*])_{(\mu)}$$ est nulle.

Pour tout $i$, on a une suite exacte $$0 \longrightarrow W_{i+1}[d \lambda^*] \longrightarrow W_i[d \lambda^*] \longrightarrow (W_i / W_{i+1})[d \lambda^*] \longrightarrow 0,$$ donc une suite exacte sur les composantes isotypiques $$ R^1 \operatorname{Ind} (W_{i+1}[d \lambda^*] )_{(\mu)} \rightarrow R^1 \operatorname{Ind} (W_i[d \lambda^*] )_{(\mu)} \rightarrow R^1 \operatorname{Ind} ((W_i / W_{i+1})[d \lambda^*] )_{(\mu)}.$$

Il suffit donc de montrer que pour tout $i$, on a $$ R^1 \operatorname{Ind} ((W_i / W_{i+1})[d \lambda^*] )_{(\mu)}=0 $$ et la proposition sera prouvée.

Supposons le contraire : selon le théorème de Borel-Weil-Bott, il existe une racine simple $\alpha$ telle que 
 $$ s _ \alpha \star \mu = d \lambda ^ * + \beta,$$ où $\beta$ est le plus grand poids de $W_i/W_{i+1}$ (c'est donc une racine de $G$ qui n'est pas une racine de $Q_\lambda $).
On en déduit $$ \mu - d \lambda ^ * =N \alpha + \beta,$$ en posant $$N :=  1 + \langle \alpha ^ \vee , \mu \rangle \geq 1$$

Comme $d<0$, le poids $N \alpha + \beta$ est dominant : la racine simple $\alpha$ et la racine positive $\beta$ sont donc données dans la liste du \S 1.3.1 (si $N=1$) ou de la proposition \ref{t3.9} (si $N\geq 2$).
Ici, on remarque de plus que :
\begin{itemize}
 \item Le poids $N \alpha + \beta$ est la somme de deux poids dominants $\mu$ et $-d \lambda^*$ non nuls.
 \item On a $\langle \alpha ^ \vee , \mu \rangle=0$ si et seulement si $N=1$.
 \item On a $\langle \beta ^ \vee , \lambda^* \rangle \neq 0$ (car $\beta$ n'est pas une racine de $Q_\lambda $).
 \end{itemize}

En consultant les deux listes, on constate immédiatement que cela est impossible (si $G$ n'est pas de type $A_1$).
\hfill $\Box$

\subsection{Cas où le groupe $G$ n'est pas simple}

On suppose dans cette partie que le groupe $G$ est de la forme $G^1 \times G^2$.
Le sous-groupe de Borel $B$ de $G$ s'écrit $B = B^1 \times B^2$~; de même pour le tore maximal $T = T^1 \times T^2$.
La représentation $V(\lambda)$ de $G$ s'écrit $V(\lambda)= V(\lambda_1) \otimes V(\lambda_2)$, où l'on suppose les poids dominants respectifs $\lambda_1$ et $\lambda_2$ de $G^1$ et $ G^2$ tous les deux non nuls.

On note $P_{\lambda_i}$ le stabilisateur dans $G^i$ de la droite des vecteurs de plus grand poids de $V(\lambda_i)$. On a $P_{\lambda} = P_{\lambda_1} \times P_{\lambda_2}$.

Notre variété de drapeaux est donc un produit $X_{\lambda} = X_{\lambda_1} \times X_{\lambda_2} $, où l'on note $X_{\lambda_i}$ la variété de drapeaux $G^i/P_{\lambda_i}$.
On note $p_i  : X_{\lambda} \longrightarrow X_{\lambda_i}$ les projections canoniques.

On a $$\Lcal_{\lambda}= p_1 ^ * \Lcal_{\lambda_1} \otimes p_2^* \Lcal_{\lambda_2}, $$
et $$\Tcal_{X_{\lambda}}= p_1 ^ * \Tcal_{X_{\lambda_1}} \oplus p_2^* \Tcal_{X_{\lambda_2}}. $$\\

On remarque que selon la proposition \ref{t3.6} $$H^1(E_\lambda, \Ocal_{E_\lambda}) = \bigoplus_{d \in \Zbb} H^1(X_\lambda,  \Lcal_{d \lambda})=0.$$ La proposition \ref{t3.3} donne donc
 $$  T^1_\lambda \cong H^1(E_\lambda, \Tcal_{E_\lambda}).$$

Cet isomorphisme est aussi conséquence de \cite{Sern} prop II.5.8 (ii). En effet, $\operatorname{dim}(C_\lambda)= \operatorname{dim}(X_{\lambda_1} )  + \operatorname{dim}(X_{\lambda_2}) +1 \geq 3$. Comme $C_\lambda$ est Cohen-Macaulay \cite{Ra}, sa profondeur en $0$ est supérieure ou égale à $3$.

\subsubsection{Cas où $G= \operatorname{SL}(2) \times \operatorname{SL}(2)$}

On note encore $H$ le groupe des matrices unipotentes triangulaires inférieures de taille $2 \times 2$.

On écrit le poids dominant de $G$ sous la forme $\lambda = (m,n)$, où l'on peut supposer les entiers $m$ et $n$ tels que $m \geq n \geq 1$.

Pour calculer $H^1(E_\lambda, \Tcal_{E_\lambda})  \cong H^1(X_\lambda, \pi _ * \Tcal_{E_\lambda})$, on va utiliser une résolution du faisceau $G$-linéarisé $\pi _ * \Tcal_{E_\lambda}$. Sa fibre en $Q_\lambda$ est le $Q_\lambda$-module $$ \bigoplus _{d \in \Zbb} \frak{g}/  \frak{g}_{v_{-\lambda}}[md,nd].$$

D'où $$H^1(E_\lambda, \Tcal_{E_\lambda}) = \bigoplus _{d \in \Zbb} R^1 \operatorname{Ind} (\frak{g}/  \frak{g}_{v_{-\lambda}}[md,nd]).$$
 
Soit $d \in \Zbb$. On remarque que $$\frak{g}_{v_{-\lambda}} = (\frak{h} \times \frak{h}) \oplus \frak{t}_{v_{-\lambda}},$$ où le stabilisateur $ \frak{t}_{v_{-\lambda}}$ de $v_{-\lambda}$ dans $ \frak{t}$ est une droite $Q_\lambda$-invariante.

On a donc une suite exacte de $Q_\lambda $-modules $$    0   \longrightarrow            \Cbb[md,nd]              \longrightarrow             \frak{g}/(\frak{h} \times \frak{h})[md,nd]                 \longrightarrow                          \frak{g}/  \frak{g}_{v_{-\lambda}}[md,nd]        \longrightarrow         0  .$$ D'où une suite exacte longue $$  ...\rightarrow   R^1 \operatorname{Ind}(      \Cbb[md,nd]  )            \rightarrow     R^1  \operatorname{Ind}(      \frak{g}/(\frak{h} \times \frak{h})[md,nd]      )           \rightarrow         R^1   \operatorname{Ind}(              \frak{g}/  \frak{g}_{v_{-\lambda}}[md,nd]  )  $$ $$    \rightarrow   R^2  \operatorname{Ind}(  \Cbb[md,nd] ) \xrightarrow{h_d}    R^2  \operatorname{Ind}(      \frak{g}/(\frak{h} \times \frak{h})[md,nd] )      \rightarrow   ...$$

\begin{prop}\label{t3.11}
\begin{enumerate}[(1)]

\item L'espace $R^1 \operatorname{Ind}(      \Cbb[md,nd]  )   $ est nul pour tout $d$.

\item L'espace $ R^1  \operatorname{Ind}(      \frak{g}/(\frak{h} \times \frak{h})[md,nd]      )  $ est nul sauf si $d=-1$ et $n=1$. Dans ce cas il vaut $V(m-2, 1)$.

\item Le noyau de $h_d$ est nul, sauf si $d=-1$ et $(m,n)=(2,2)$ et si $d=-2$ et $(m,n)=(1,1)$. Dans ces deux cas (qui correspondent au cas (H9) du théorème 1.1) il vaut $V(0, 0) $.

\end{enumerate}

\end{prop}

\noindent {\bf Preuve.} 

(1) Cela découle immédiatement du théorème de Borel-Weil-Bott (ou simplement de la cohomologie des faisceaux inversibles sur $\Pbb ^1 \times \Pbb ^1$).

(2) On remarque qu'on a l'isomorphisme de $Q_\lambda $-modules
$$\frak{g}/(\frak{h} \times \frak{h})[md,nd] \simeq  V(0,1)[md,nd+1] \oplus V(1,0)[md+1,nd],$$ avec par exemple $$R^1  \operatorname{Ind}( V(0,1)[md,nd+1] ) = V(0,1)  \otimes R^1  \operatorname{Ind} (\Cbb[md,nd+1]).$$ 

Selon le théorème de Borel-Weil-Bott, pour que l'espace $R^1  \operatorname{Ind} (\Cbb[md,nd+1])$ soit non nul, il faut que les entiers $md$ et $nd+1$ soient l'un positif, l'autre strictement négatif. Comme $md$ et $nd$ sont de même signe, on a nécessairement $md<0$ et $nd+1 \geq 0$ (donc $d=-1$ et $n=1$).
Dans ce cas, $$ R^1  \operatorname{Ind}(V(0,1)[md,nd+1])=  V(0,1)  \otimes V(m-2,0) = V(m-2,1).$$
Il faut donc que l'on ait $m\geq 2$.

De la même façon, on obtient que $R^1  \operatorname{Ind}( V(1,0)[md+1,nd])$ est toujours nul, car on a supposé $m \geq n$.

(3) On va en fait déterminer le conoyau de l'application transposée$ ~^t h_d$. Selon le théorème de dualité de Serre (\cite{H} III.7), on a un isomorphisme fonctoriel                                                              $$R^2  \operatorname{Ind}(W)^* \cong \operatorname{Ind}(W^*[-2,-2])$$ pour tout $Q_\lambda$-module $W$ (car la fibre en $Q_\lambda$ du faisceau anticanonique de $X_\lambda$ est $\Cbb[-2,-2]$).\\

Supposons que le conoyau de l'application $$\operatorname{Ind}([{\frak{g}/(\frak{h} \times \frak{h})] }^*[-md-2,-nd-2]) \xrightarrow{ ~^t h_d~} \operatorname{Ind}(\Cbb[-md-2,-nd-2]).$$ est non nul.

Pour que l'espace d'arrivée de $~^t h_d$ soit non nul, il faut que $M :=-md-2$ et $N :=-nd-2$ soient positifs. Cet espace est alors le module simple $V(M,N)$, et il faut que l'application $ ~^t h_d$ soit nulle.

Comme dans la démonstration de (2), on remarque que l'espace de départ de $ ~^t h_d$ est une somme directe : $$\operatorname{Ind}([{\frak{g}/(\frak{h} \times \frak{h})] }^*[-md-2,-nd-2]) \simeq V(0,1)[M,N-1] \oplus V(1,0)[M-1,N].$$

Il faut donc que les deux composantes de $~^t h_d$ soient nulles.
La première composante est l'image par $\operatorname{Ind}$ du morphisme $j$ de la suite exacte courte suivante
$$ 0 \xrightarrow{~~~~} \Cbb[M,N-2] \xrightarrow{~~~~} V(0,1)[M,N-1] \xrightarrow{~~j~} \Cbb[M,N] \xrightarrow{~~~~} 0.$$

Le conoyau de la première composante se plonge donc dans l'espace $R^1 \operatorname{Ind}(\Cbb[M,N-2])$. Selon le théorème de Borel-Weil-Bott, pour que ce dernier espace soit non nul, il faut $N-2 \leq -2$, donc $N=0$. On montre de même, à l'aide de la seconde composante de $~^t h_d$, que $M=0$.

Enfin, dans le cas où $(M,N)=(0,0)$, le conoyau de $ 0 \xrightarrow{ ~^t h_d~} V(0)$ est bien $V(0)$.
\hfill $\Box$

\subsubsection{Autres cas}

Lorsque les actions de $G^1$ et $G^2$ se factorisent par $\operatorname{SL}(2)$ : $$G^i \longrightarrow \operatorname{SL}(2) \longrightarrow \operatorname{GL}(V(\lambda_i)),$$ on est dans la situation du \S3.3.1~; on suppose donc que l'action de $G^2$ ne se factorise pas par $\operatorname{SL}(2)$.

\begin{prop} \label{t3.12}

\begin{enumerate}[(1)]

\item L'espace $H^2(X_\lambda , \Lcal_{d \lambda} )$ est nul pour tout entier $d$.

\item On a donc un isomorphisme $ T^1_\lambda \cong \bigoplus_{d \in \Zbb} H^1(X_\lambda , \Lcal_{d \lambda} \otimes \Tcal_{X_\lambda})$.

\end{enumerate}

\end{prop}

\noindent {\bf Preuve.}

(1) On rappelle que si $d \geq 0$, tous les groupes de cohomologie de $\Lcal_{d \lambda}$ sont nuls sauf en degré $0$, et si $d<0$, le groupe de cohomologie de degré $0$ est nul.
On peut donc supposer $d<0$, et on a, selon la formule de Künneth (\cite{Da} p32) $$H^2(X_\lambda , \Lcal_{d \lambda} ) = H^1(X_{\lambda_1} , \Lcal_{d \lambda_1}) \otimes H^1(X_{\lambda_2} , \Lcal_{d \lambda_2}). $$ Or selon la proposition \ref{t3.6} et l'hypothèse faite sur $\lambda_2$, on a $H^1(X_{\lambda_2} , \Lcal_{d \lambda_2})=0$, d'où le résultat.

(2) Selon le point (1) et la proposition \ref{t3.6}, la suite exacte (\ref{ty8}) s'écrit

$$ 0 \longrightarrow  H^1(E_\lambda,\Tcal_{E_\lambda})  \longrightarrow  \bigoplus_{d \in \Zbb} H^1(X_\lambda , \Lcal_{d \lambda} \otimes \Tcal_{X_\lambda}) \longrightarrow 0,$$ d'où le résultat.
\hfill $\Box$

La proposition suivante achève donc la démonstration du théorème 3.1 :

\begin{prop} \label{t3.13}  L'espace $H^1(X_\lambda , \Lcal_{d \lambda} \otimes \Tcal_{X_\lambda})$ est nul, sauf dans les cas suivants (à factorisation près) :

\begin{itemize}

 \item Si $G=\operatorname{SL}(2) \times \operatorname{SL}(n)$ et $d=-1$ et $\lambda = ( m  , \omega _ 1)$, alors il vaut 
$V_{\operatorname{SL}(2)}(m-2) \otimes V(\omega_1)$.
 
 \item Si $G=\operatorname{SL}(2) \times \operatorname{SL}(n)$ et $d=-1$ et $\lambda = ( m  , \omega _ n)$, alors il vaut 
$V_{\operatorname{SL}(2)}(m-2) \otimes V(\omega_n)$.
 
\item Si $G=\operatorname{SL}(2) \times \operatorname{Sp}(2n)$ et $d=-1$ et $\lambda = ( m  , \omega _ 1)$, alors il vaut 
$V_{\operatorname{SL}(2)}(m-2) \otimes V(\omega_1)$.

 \end{itemize}

\end{prop}

\noindent {\bf Preuve.} 
Selon la proposition \ref{t3.7}, on peut supposer $d<0$.
On a $$\Lcal_{d \lambda} \otimes \Tcal_{X_\lambda} = (p_1^* (\Lcal_{d \lambda_1} \otimes \Tcal_{X_{\lambda_1}}   ) \otimes p_2^* (\Lcal_{d \lambda_2}  ) ) \oplus ( p_1^* ( \Lcal_{d \lambda_1}   ) \otimes p_2^* (\Lcal_{d \lambda_2} \otimes \Tcal_{X_{\lambda_2}}  )  ).$$

Donc, selon la formule de Künneth, $$H^1(X_\lambda , \Lcal_{d \lambda} \otimes \Tcal_{X_\lambda})= H^1(X_{\lambda_1}, \Lcal_{d \lambda_1} ) \otimes H^0(X_{\lambda_2},  \Lcal_{d \lambda_2} \otimes \Tcal_{X_{\lambda_2}} ).$$ En effet, comme $d<0$, on a $$H^0(X_{\lambda_1}, \Lcal_{d \lambda_1} )= H^0(X_{\lambda_2}, \Lcal_{d \lambda_2} )=0$$ et selon la proposition \ref{t3.6} (grâce à l'hypothèse faite sur $\lambda_2$), on a $$H^1(X_{\lambda_2}, \Lcal_{d \lambda_2} ) =0.$$ 

Si l'action de $G^1$ sur $V(\lambda_1)$ ne se factorise pas par $\operatorname{SL}(2)$, on aura également $H^1(X_{\lambda_1}, \Lcal_{d \lambda_1} ) =0$. On peut donc supposer que $G^1 = \operatorname{SL}(2)$, de sorte que  $H^1(X_{\lambda_1}, \Lcal_{d \lambda_1} ) $ est non nul, et vaut $V(-d \lambda_1 -2)$ (on considère désormais $\lambda_1$ comme un entier).\\

Il reste à calculer l'espace des sections globales $H^0(X_{\lambda_2},  \Lcal_{d \lambda_2} \otimes \Tcal_{X_{\lambda_2}} )$.

Remarquons tout d'abord que sa partie $G$-invariante est nulle : en effet, selon l'isomorphisme (où $\Cbb$ est la droite munie de l'action triviale de $G$) $$\operatorname{Hom}^G(\Cbb, \operatorname{Ind}( \frak{g}^2/ {\frak{q}} _{\lambda_2}[d \lambda_2^* ])) \cong \operatorname{Hom}^{Q_{\lambda_2}}(\Cbb, \frak{g}^2/  {\frak{q}}_{\lambda_2}[d \lambda_2^* ]),$$ la partie $G$-invariante est isomorphe à l'espace des $Q_{\lambda_2}$-invariants suivant :
$$ (\frak{g}^2/  {\frak{q}}_{\lambda_2}[d \lambda_2^* ])^{Q_{\lambda_2}}.$$

Supposons par l'absurde ce dernier espace non nul. Ses éléments sont en particulier de poids nul pour le tore $T$ : ils admettent donc un représentant dans $\frak{g}^2$ dont le poids est une racine positive $\beta$ de $\frak{g}^2$ telle que $$\beta + d \lambda_2^* =0.$$ On en déduit d'une part que l'on peut supposer que le groupe $G^2$ est simple (car son action sur $V(\lambda_2)$ se factorise par celle d'un groupe simple), et d'autre part que la racine positive $\beta$ est une racine dominante.
Or on sait qu'alors le sous-$ {\frak{q}}_{\lambda_2}$-module de $\frak{g}^2$ engendré par  $\frak{g}^2_\beta$ contient $\frak{g}^2_\alpha$ pour toute racine simple $\alpha$ de $\frak{g}^2$.

Ainsi, toute racine simple de $\frak{g}^2$ est une racine de $ {\frak{q}}_{\lambda_2}$, et $ {\frak{q}}_{\lambda_2}=\frak{g}^2$ : une contradiction.\\

Déterminons maintenant les autres composantes isotypiques de l'espace des sections globales : soit $\mu$ un poids dominant de $G^2$ non nul. Supposons que $H^0(X_{\lambda_2},  \Lcal_{d \lambda_2} \otimes \Tcal_{X_{\lambda_2}} )_{(\mu)} $ est non nulle.

Comme dans la démonstration de la proposition \ref{t3.10}, appliquons le théorème de Borel-Weil-Bott à l'aide d'une suite de Jordan-Hölder du $ Q_{\lambda_2}$-module $ \frak{g}^2/\frak{q}_{\lambda_2}$.

On obtient qu'il existe une racine positive $\beta$ de $G^2$ telle que $$\beta + d \lambda^*_2 = \mu.$$

Comme $d<0$, pour que $\beta + d \lambda^*_2$ soit dominant, il faut que l'action de $G^2$ sur $V(\lambda_2)$ se factorise par un groupe simple~; on suppose donc que $G^2$ est un groupe simple.

On remarque que $\beta = \mu - d \lambda^*_2$ est un poids dominant, non fondamental (car $\mu$ et $- d \lambda^*_2$ sont tous les deux non nuls).

Or les seules racines dominantes des systèmes de racines irréductibles qui ne sont pas des poids fondamentaux sont   
\begin{itemize}
 \item la plus grande racine $\omega_1 + \omega_n$ des systèmes de racines de type $A_n$ ($n \geq 1$).
 \item la plus grande racine $2 \omega_1$ des systèmes de racines de type $C_n$ ($n \geq 2$).  
 \end{itemize}

On a donc $d=-1$, et l'on est dans l'une des deux situations suivantes (à revêtement fini de $G^2$ près) :
\begin{itemize}
 \item On a $G^2 = \operatorname{SL}(n)$ ($n\geq 2$), et $\mu = \lambda_2 = \omega_1$ ou $\mu = \lambda_2 = \omega_n$.
 \item On a $G^2 = \operatorname{Sp}(2n)$ ($n\geq 2$), et $\mu = \lambda_2 = \omega_1$.
 \end{itemize}

Il ne reste plus qu'à vérifier que dans ces deux cas, l'espace $ T^1_\lambda$ est celui annoncé (la seule autre possibilité est qu'il soit nul).

En utilisant les notations analogues à celles de la démonstration de la proposition \ref{t3.10} (en remplaçant $G$ par $G^2$), on a pour tout $i$ une suite exacte 
$$0 \rightarrow  \operatorname{Ind} (W_{i+1}[d \lambda_2^*] )_{(\mu)}  \rightarrow   \operatorname{Ind} (W_i[d \lambda_2^*] )_{(\mu)}       \rightarrow  \operatorname{Ind} (W_i / W_{i+1}[d \lambda_2^*] )_{(\mu)}   \rightarrow  R^1 \operatorname{Ind} (W_{i+1}[d \lambda_2^*] )_{(\mu)} .$$

Selon cette même démonstration, on a $$ R^1 \operatorname{Ind} (W_{i+1}[d \lambda_2^*] )_{(\mu)}=0$$ pour tout $i$.
On en conclut facilement que $$H^0(X_{\lambda_2},  \Lcal_{d \lambda_2} \otimes \Tcal_{X_{\lambda_2}} )= \operatorname{Ind} (W_{0}[d \lambda_2^*] )= V(\lambda_2)$$ dans les deux situations.
\hfill $\Box$

\subsection{Démonstration de la proposition 3.3}

On va d'abord démontrer que la déformation $\frak{V}$ de la proposition 3.3 est
plate (proposition \ref{defplate}). Elle est donc déduite de la déformation verselle de $C_{mn}$ par un
changement de base. On montrera ensuite (proposition \ref{defnontriv}) que la différentielle du changement de
base est un isomorphisme (de l'espace tangent à $\Spec \Cbb[[\mathbf{t}]]$ vers l'espace
des déformations infinitésimales de $C_{mn}$), ce qui
démontrera la proposition 3.3.\\

Soit $N$ un entier supérieur ou égal à $2$.
On munit l'ensemble des mon\^omes de l'anneau de polynomes
$\Cbb[x_1,...,x_N,z]$ de l'ordre lexicographique : si $m_1=x_1^{\alpha_1}...x_N^{\alpha_N}z^{\alpha_{N+1}}$ et $m_2=
x_1^{\beta_1}...x_N^{\beta_N}z^{\beta_{N+1}}   $ sont deux  mon\^omes distincts,
alors $m_1<m_2$ si et seulement si $\alpha_i<\beta_i$ pour le plus petit
indice $i$ tel que $\alpha_i \neq \beta_i$.


On renvoie à \cite{Eis} p 332 pour la définition d'une base de Gr\"obner. Le
lemme suivant est une application immédiate du critère de Buchberger (\cite{Eis} p 338):

\begin{lem}
Soient $s_1, ... , s_N$ des nombres complexes.
Notons $g_{ij}$ le mineur $2 \times 2$ de la matrice 

$$\displaystyle{\left(\begin{array}{cccc} x_{0}  & x_{1}  &  \cdot \cdot \cdot &
      x_{N-1}  \\ x_{1}-s_1z & x_{2}-s_2z &  \cdot \cdot \cdot &  x_{N}-s_Nz
      \end{array}\right) }$$
obtenu en prenant les colonnes $i$ et $j$.
La famille $(g_{ij})_{i<j}$ est une base de Gr\"obner de l'idéal qu'elle engendre.
\end{lem}

On reprend maintenant les notations de la proposition 3.3.

\begin{prop}\label{defplate}
Notons  $\frak{X}$ la famille de sous schémas fermés $$\xymatrix{ \frak{X}
  \ar[rd] &   \subseteq &  \Spec\Cbb[\mathbf{t},\mathbf{x}]=\Abb ^ {(m-1)(n+1)} \times \Abb ^
  {(m+1)(n+1)} \ar[ld] ^ {pr_1} \\ & ~~~~~\Spec\Cbb[\mathbf{t}]=\Abb ^ {(m-1)(n+1)}   } $$ définie par
  les mineurs $2  \times  2$ de (\ref{verse}).
La famille  $\frak{X}$ est plate au dessus de $ \Abb ^ {(m-1)(n+1)} $.\\
Par conséquent la famille $\frak{V}$  est plate.

\end{prop}

\noindent {\bf Preuve.} 
On va montrer que l'adhérence schématique  $\overline{ \frak{X}   }$  de $\frak{X}$
dans $ \Abb ^ {(m-1)(n+1)} \times \Pbb ^ {(m+1)(n+1)}$ est une famille plate
sur  $ \Abb ^ {(m-1)(n+1)} $, ce qui donnera le résultat.

Pour cela, il suffit (\cite{H} thm 9.9 p 261) de montrer que le polyn\^ome de
Hilbert de la fibre $\overline{  \frak{X}  }_s$ de  $\overline{ \frak{X}   }$  au dessus de $s=(s_{ij})_{i,j}$
(vue comme un sous-schéma fermé de $\Pbb ^ {(m+1)(n+1)}$) est indépendant de
  $s$.

L'idéal homogène $I_s$ de $\overline{\frak{X}}_s$ est engendré par les mineurs
$2  \times  2$ de la matrice
$$ \left(\begin{array}{cccc|c|cccc} x_{00} &
     ... &  x_{0~m-2} &
      x_{0~m-1} &   \cdot \cdot \cdot &  x_{n0} &... & x_{n~m-2} & x_{n~m-1} \\ x_{01}-s_{01}z & ...  &  x_{0~m-1} - s_{0~m-1}z &
      x_{0~m} &  \cdot \cdot \cdot &  x_{n1} -s_{n1}z &... & x_{n~m-1} -
     s_{n~m-1}z & x_{n~m}  \end{array}\right)         $$ où $z$ est une
     indéterminée supplémentaire. 
Selon le lemme précédent, ceux-ci forment une base de Gr\"obner pour
     l'ordre lexicographique tel que $x_{ij}<x_{kl}$ si $i<k$ ou si $i=k$ et
     $j<l$, et tel que $z$ soit la plus petite indéterminée. 

On remarque que les termes initiaux de ces mineurs $2  \times  2$ ne dépendent
pas de $s$. Selon \cite{Eis} Thm 15.3 p 329, le $\Cbb$-espace vectoriel
quotient $\Cbb[\mathbf{x}, z]/I_s$ admet donc pour base les classes
modulo $I_s$ d'une famille de mon\^omes indépendante de $s$.
La dimension des composantes homogènes de $\Cbb[\mathbf{x}, z]/I_s$ est
donc indépendante de $s$, et le polyn\^ome de Hilbert de $\overline{  \frak{X}
}_s$ aussi, ce qui démontre la proposition.
\hfill $\Box$\\

Remarquons que l'espace tangent à $\Spec \Cbb[[\mathbf{t}]]$ est de m\^eme
dimension  que $T^1(C_{mn})$. Pour montrer que la différentielle du morphisme
correspondant à $\frak{V}$
(de $\Spec \Cbb[[\mathbf{t}]]$ vers la déformation verselle de $C_{mn}$) est un
isomorphisme entre les espaces tangents, il
suffit donc de montrer qu'elle est injective, c'est-à-dire que tout vecteur tangent non nul à $\Spec \Cbb[[\mathbf{t}]]$ correspond à une
déformation infinitésimale de $C_{mn}$ non triviale.

Notons $\epsilon$ la classe de $y$ dans l'algèbre $\Cbb[y]/\langle y^2 \rangle$.

\begin{prop} \label{defnontriv}
Soit un vecteur tangent à $\Spec \Cbb[[\mathbf{t}]]$, c'est-à-dire un
morphisme de $\Spec \Cbb[\epsilon]$ vers $\Spec \Cbb[[\mathbf{t}]]$.
On note $s_{ij}$ les nombres complexes tels que le morphisme correspondant\\  $  \Cbb[[\mathbf{t}]] \longrightarrow
\Cbb[\epsilon]  $ envoie $t_{ij}$ vers $\epsilon s_{ij}$.

La déformation induite par $\frak{V}$ sur $\Spec \Cbb[\epsilon]$ est le
sous-schéma fermé de $\Spec \Cbb[\mathbf{x}][\epsilon]$ défini par les mineurs
$2  \times  2$ de la matrice $$ \left(\begin{array}{cccc|c|cccc} x_{00} &
     ... &  x_{0~m-2} &
      x_{0~m-1} &   \cdot \cdot \cdot &  x_{n0} &... & x_{n~m-2} & x_{n~m-1}
      \\ x_{01}-\epsilon s_{01} & ...  &  x_{0~m-1} -\epsilon s_{0~m-1} &
      x_{0~m} &  \cdot \cdot \cdot &  x_{n1} -\epsilon s_{n1} &... & x_{n~m-1}
      -\epsilon s_{n~m-1} & x_{n~m}  \end{array}\right)  .$$ Elle est triviale
  si et seulement si les $s_{ij}$ sont tous nuls. 
\end{prop}

\noindent {\bf Preuve.} On va appliquer \cite{H} ex 9.8 p 267.
Notons $J$ l'idéal du sous-schéma fermé $C_{mn}$ de $\Abb ^ {(m+1)(n+1)}$, et
$A:=\Spec \Cbb[\mathbf{x}]/J$ son algèbre affine. La déformation
infinitésimale de la proposition correspond au morphisme de $A$-modules $\phi
: J/J^2
\longrightarrow A $ tel que
$$ \phi \left(\begin{array}{c|cc|c} &x_{ij} &
       x_{kl}& \\ &x_{i~j+1} &  x_{k~l+1}& \end{array}\right)=-\begin{array}{|cc|} x_{ij} &
       x_{kl} \\ s_{i~j+1} &  s_{k~l+1} \end{array} ~ .$$

Elle est triviale si et seulement si le morphisme $\phi $ est induit par un
champ de vecteurs $$\sum_{i,j}h_{ij} \frac{\partial}{\partial x_{ij}}, $$
c'est-à-dire s'il existe des des éléments $h_{ij}$ de $A$ tels que
$$ \phi \left(\begin{array}{c|cc|c} &x_{ij} &
       x_{kl}& \\ &x_{i~j+1} &  x_{k~l+1}& \end{array}\right)=\begin{array}{|cc|} x_{ij} &
       x_{kl} \\ h_{i~j+1} &  h_{k~l+1} \end{array}+ \begin{array}{|cc|} h_{ij} &
       h_{kl} \\ x_{i~j+1} &  x_{k~l+1} \end{array} ~ .$$

En prenant alors $i$ et $k$ distincts (ce qui est possible car $n$ est
non nul), on obtient que la composante homogène de degré $1$ de $h_{ij}$ est nulle dès que $j<m$, puis que les $s_{ij}$
sont tous nuls, d'où le résultat. 
\hfill $\Box$

\section*{Appendice}

On note $\operatorname{Hilb}(\Pbb(V(\lambda)))$ le schéma de Hilbert (construit dans \cite{Gr}) des sous-schémas fermés de $\Pbb(V(\lambda))$. Le sous-schéma $\operatorname{Hilb}^G(\Pbb(V(\lambda)))$  des points fixes de $G$ dans $\operatorname{Hilb}(\Pbb(V(\lambda)))$ paramètre les sous-schémas fermés de $\Pbb(V(\lambda))$ qui sont stables par $G$.

On répond dans cet appendice à la question naturelle suivante : quelles déformations locales du cône des vecteurs primitifs $C(\lambda)$ peut-on obtenir à l'aide de $\operatorname{Hilb}(\Pbb(V(\lambda)))$?

Comme $C(\lambda)$ est le cône affine dans $V(\lambda)$ au dessus de la variété de drapeaux $$X_\lambda  := G/P_\lambda \subseteq \Pbb(V(\lambda)),$$ on peut être tenté de déformer $X_\lambda$ dans $\Pbb(V(\lambda))$ à l'aide du schéma de Hilbert pour en déduire naturellement une déformation de $C(\lambda)$.

La proposition suivante montre que l'on n'obtient ainsi que des déformations triviales, c'est-à-dire provenant de l'action du groupe $\operatorname{GL}(V(\lambda))$ des automorphismes d'espace vectoriel de $V(\lambda)$.

On note $z$ le point de $\operatorname{Hilb}(\Pbb(V(\lambda)))$ correspondant à $X_\lambda$.
Le groupe $\operatorname{GL}(V(\lambda))$ agit naturellement sur $\operatorname{Hilb}(\Pbb(V(\lambda)))$.

\begin{prop} L'orbite $\operatorname{GL}(V(\lambda)) \cdot z$ est ouverte dans $\operatorname{Hilb}(\Pbb(V(\lambda)))$.
\end{prop}

\noindent {\bf Preuve.} 

Il suffit de montrer que l'espace tangent à l'orbite est égal à l'espace tangent au schéma de Hilbert :
$$T_z (\operatorname{GL}(V(\lambda)) \cdot z )     =    T_z \operatorname{Hilb}(\Pbb(V(\lambda))),$$
c'est-à-dire que l'application $$\frak{gl}(V(\lambda)) \xrightarrow{~\phi~} T_z \operatorname{Hilb}(\Pbb(V(\lambda)))$$ obtenue en différentiant l'application naturelle  $$\begin{array}{clc} \operatorname{GL}(V(\lambda)) & \longrightarrow&  \operatorname{Hilb}(\Pbb(V(\lambda))),\\ u &\longmapsto & u \cdot z
\end{array}$$  est surjective.

Notons $\Ncal_{X_\lambda}$ le faisceau normal à $X_\lambda$ dans $\Pbb(V(\lambda))$.

Il est donné par la suite exacte courte de $\Ocal_{X_\lambda}$-modules : \begin{equation} \label{sec1} 0 \longrightarrow \Tcal_{X_\lambda} \longrightarrow \Tcal_{\Pbb(V(\lambda))}|_{X_\lambda}   \longrightarrow  \Ncal_{X_\lambda}   \longrightarrow  0.\end{equation}

L'espace tangent au schéma de Hilbert en $z$ est canoniquement isomorphe à l'espace des sections globales de $\Ncal_{X_\lambda}$ : $$ T_z \operatorname{Hilb}(\Pbb(V(\lambda))) \cong \operatorname{H}^0(X_\lambda,\Ncal_{X_\lambda}).$$ 

On sait que l'espace $$\operatorname{H}^1(X_\lambda, \Tcal_{X_\lambda})$$ est nul (cela résulte par exemple de la proposition \ref{t3.7}).

En utilisant la suite exacte courte (\ref{sec1}), on en déduit que l'application canonique $$\operatorname{H}^0(X_\lambda,\Tcal_{\Pbb(V(\lambda))}|_{X_\lambda}) \xrightarrow{~\phi_1~}  \operatorname{H}^0(X_\lambda,\Ncal_{X_\lambda})$$ est surjective.

On utilise ensuite la suite exacte de $\Ocal_{\Pbb(V(\lambda))}$-modules (\cite{H} Example II.8.20.1) :
$$0    \longrightarrow  \Ocal_{\Pbb(V(\lambda))}   \longrightarrow  \Ocal_{\Pbb(V(\lambda))}(1)\otimes_\Cbb V(\lambda)   \longrightarrow  \Tcal_{\Pbb(V(\lambda))}    \longrightarrow  0.$$

Comme les termes de cette suite sont des faisceaux localement libres, on obtient encore une suite exacte si on la restreint à $X_\lambda$ :
$$0    \longrightarrow  \Ocal_{X_\lambda}   \longrightarrow  \Lcal_\lambda \otimes_\Cbb V(\lambda)   \longrightarrow  \Tcal_{\Pbb(V(\lambda))}|_{X_\lambda}    \longrightarrow  0.$$

L'application associée $$V(\lambda)^* \otimes_\Cbb V(\lambda) = \operatorname{H}^0(X_\lambda,\Lcal_\lambda) \otimes_\Cbb V(\lambda) \xrightarrow{~\phi_2~} \operatorname{H}^0(X_\lambda,\Tcal_{\Pbb(V(\lambda))}|_{X_\lambda})$$ est surjective, car l'espace $\operatorname{H}^1(X_\lambda,\Ocal_{X_\lambda})$ est nul (cela découle par exemple du théorème de Borel-Weil-Bott).

On remarque enfin que l'application $\phi$ s'identifie à la composée $\phi_1 \circ \phi_2$ : elle est donc surjective, d'où la proposition.
\hfill $\Box$\\

On en déduit immédiatement le corollaire suivant, qui montre que l'on n'obtient aucune déformation $G$-invariante à l'aide de $\operatorname{Hilb}(\Pbb(V(\lambda)))$.

\begin{cor} Le point $z$ est un point isolé réduit de $\operatorname{Hilb}^G(\Pbb(V(\lambda)))$.
\end{cor}

\noindent {\bf Preuve.} 
Il suffit de montrer que l'espace tangent à $\operatorname{Hilb}^G(\Pbb(V(\lambda)))$ en $z$ est nul.
On reprend les notations de la démonstration de la proposition précédente.
On a vu que l'application $$\frak{gl}(V(\lambda)) \xrightarrow{~\phi~} T_z \operatorname{Hilb}(\Pbb(V(\lambda)))$$ est surjective. Comme le groupe $G$ est réductif, on en déduit une surjection sur les espaces des $G$-invariants : $$\frak{gl}(V(\lambda))^G \xrightarrow{~\phi~} T_z \operatorname{Hilb}^G(\Pbb(V(\lambda))).$$

Or l'application $\phi$ est nulle sur l'espace $ \frak{gl}(V(\lambda))^G \cong \Cbb$ (qui est le centre de $ \frak{gl}(V(\lambda))$), d'où le corollaire.
\hfill $\Box$\\

Pinkham utilise dans \cite{Pink} \S4-5 de manière plus concluante le
schéma de Hilbert $\operatorname{Hilb}(\overline{V(\lambda)})$ des
sous-schémas fermés de l'espace projectif $\overline{V(\lambda)}$
obtenu en complétant $V(\lambda)$. Il étudie ainsi plus généralement
les déformations des cônes affines sur les variétés projectives
lisses, et montre sous certaines hypothèses (qui sont vérifiées dans notre situation) qu'on les obtient toutes
en déformant à l'aide de
$\operatorname{Hilb}(\overline{V(\lambda)})$ le complété du cône.

\end{document}